\documentclass[12pt]{amsart}
\usepackage{amsfonts,amssymb}
\usepackage{graphicx}
\usepackage{amscd}
\usepackage{eucal}
\input epsf
\newtheorem{thm}{\bf Theorem}[section]
\newtheorem{cor}[thm]{\bf Corollary}
\newtheorem{lem}[thm]{\bf Lemma}
\newtheorem{prop}[thm]{\bf Proposition}

\theoremstyle{definition}

\theoremstyle{remark}

\def\R{\mathbb{R}}

\def\Z{\mathbb{Z}}

\def\SL{{\rm SL}}
\def\PSL{{\rm PSL}}
\def\<{\langle }
\def\>{\rangle }

\def\f{\mathbf{f}}
\def\g{\mathbf{g}}

\begin{document}
\setlength{\baselineskip}{18pt}
\title[ ]
{ {Symmetries  of  quadratic forms classes and  of quadratic surds
continued fractions.
\\ Part I: A Poincar\'e model for the de Sitter world } }

\author[F.~Aicardi]{Francesca Aicardi}
\address[F.~Aicardi]{SISTIANA 56 PR (Trieste), Italy}
\email{aicardi@sissa.it }
%
%
%
%

\begin{abstract}
The problem of the classification  of the  indefinite  binary
quadratic forms  with integer coefficients is solved introducing a
special partition of the de Sitter world, where the coefficients
of the forms lie, into separate domains. Every class of indefinite
forms, under the action of the  special linear group acting on the
integer plane lattice, has a finite and well defined number of
representatives inside each one of such domains. This property
belongs exclusively to rational points on the one-sheeted
hyperboloid.

In the second part we will show how to obtain the symmetry type of
a class as well as its number of points  in all domains from a
sole representative of that class.
\end{abstract}
\maketitle

\section*{Introduction}\label{intro}

In this paper by {\it form} we mean  a binary quadratic form:
\begin{equation}\label{form}
 f= mx^2 + n y ^2 +k xy \end{equation}
where $m,n$ and $k$ are integers and $(x,y)$ runs on the integers
plane lattice.

The integer number
\[  \Delta=k^2-4mn \]
is called the {\it discriminant} of the form (\ref{form}).

Following \cite{Ar1}, we  call   a form {\it elliptic} if $\Delta<0$, {\it
hyperbolic} if $\Delta>0$ and  {\it parabolic} if $\Delta=0$.

According to the usual terminology, the
elliptic forms are said {\sl definite} and, and {\sl indefinite} the hyperbolic ones.

The problem of classifying and counting the orbits of binary
quadratic forms under the action of ${\rm SL}(2,\mathbb{Z})$ on
the $(x,y)$-plane dates back to Gauss and Lagrange
(\cite{Gau},\cite{Lag}) and was recently re-proposed by Arnold in
\cite{Ar1}.

The description  of the orbits of the positive definite forms by
means of the action of the modular group on the Poincar\'e model
of the Lobachevsky disc is well known: in this model, there is a
special tiling of the disc such that every tile is in one-to-one
correspondence with an element of the group, namely,  the element
that sends  the fundamental domain  to it. The upper sheet of the
two-sheeted hyperboloid where the coefficients of positive
definite forms lie is represented by the Lobachevsky disc, so that
every class of forms has one and only one representative in each
domain.

The complement to the plane of the Lobachevsky disc, representing the
hyperbolic forms, is not tiled by the same net of lines (for instance, the straight lines
of the Klein model,  separating the domains of the Lobachevsky disc)  into
domains of finite  area.

In this article I show, however,  that it is possible to introduce
a special partition into  separate domains of the one sheeted
hyperboloid where the coefficients of the forms lie: in each  of
such domains every orbit has a finite -- well defined -- number of
points\footnote{This is a surprising fact. Indeed, the orbit of a
generic point (i.e., with irrational coordinates) on the de Sitter
world is dense, as Arnold proved \cite{Ar2}  (our results imply
only that in each domain the number of points of such an orbit is
unbounded).}.

The situation is, however, intrinsically different from that of
the Lobachevsky disc, where  there is no distinguished fundamental
domain, i.e., all domains of the partition are equivalent.   In  our
partition of  the  de Sitter  world  {\sl there are two special domains},
that we call {\sl fundamental}.  An  $\SL(2,\Z)$ change  of the system of coordinates (in the plane $(x,y)$
of the forms,
and, consequently, on the hyperboloid)  changes the shape of {\sl a finite subset of the partition's tiles} (including  the
shape of the fundamental domains), but preserves
all the peculiar properties of the  tiling:

1) The  complement to the fundamental domains of the hyperboloid  is  separated,
by  the fundamental domains, into four regions, two of which (called {\sl upper regions}) are bounded from the
circle at $+ \infty$ of the hyperboloid and the other two regions  ({\sl down regions}) are bounded from the
circle at $-\infty$  (note that the circles at infinite are invariant under the action of the
group).

2) Each one of the two upper regions and of the down regions are partitioned into  a  countable set of domains which are in one-to-one correspondence  with all  elements  of the  semigroup of $\SL(2,\Z)$ generated by
$A=(^1_0 \ ^1 _1)$ and $B=(^1_1 \ ^0 _1)$.

3) Every orbit  has  the same  finite number of points, say $N_u$, in each domain of the upper regions,  and the same finite number of points, say $N_d$, in each domain of the down regions.
The fundamental domains contain   $N=N_u+N_d$ of that orbit.

To  understand  this unusual  situation (where  the  partition changes without changing the number
of integer points in the corresponding domains)    we give an example  where the reader may verify the properties
1 and 3  above (see also  Figure \ref{capp12} to see more points of the orbits, and more  domains).

 \centerline{\epsfbox{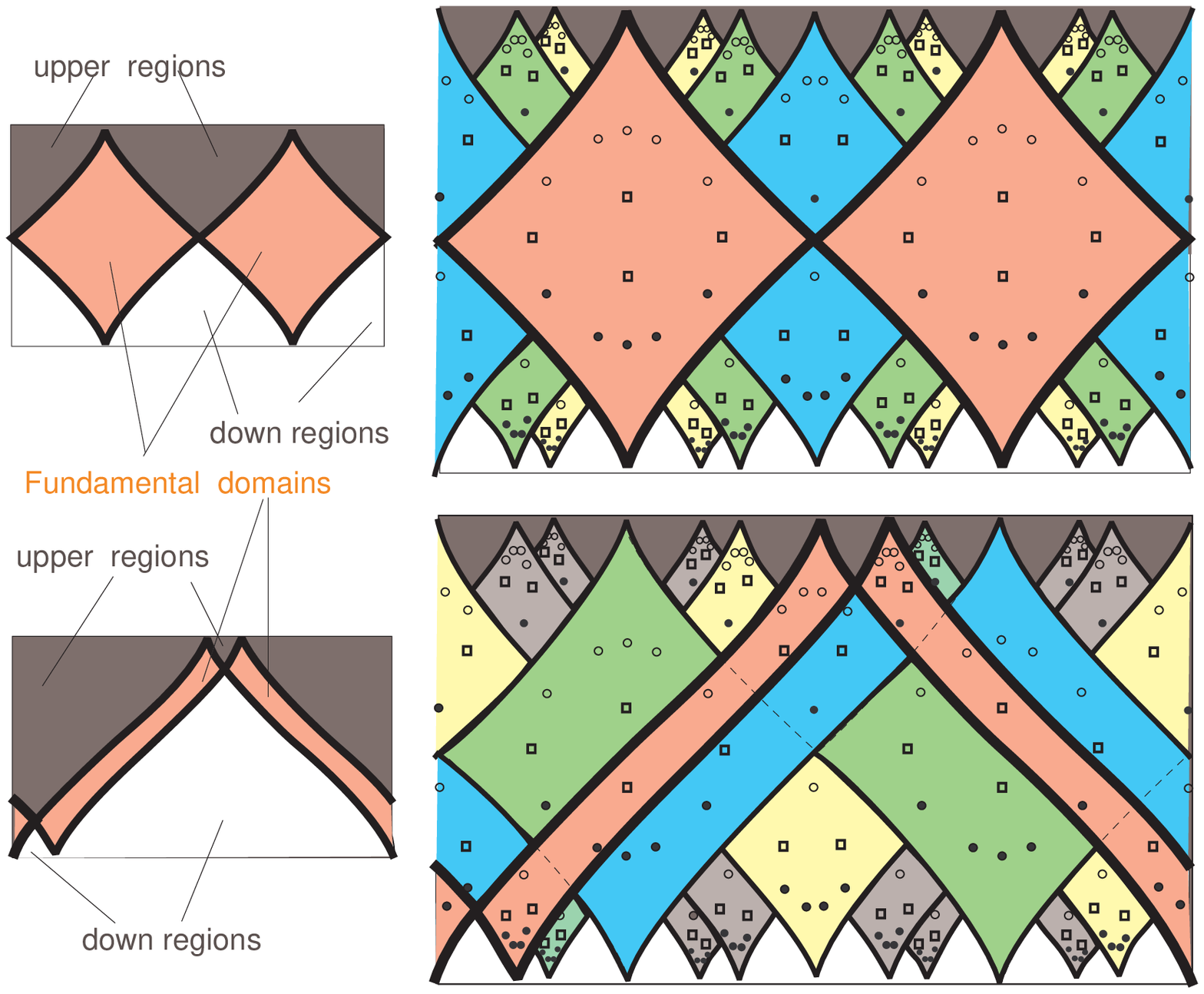}}

This  figure shows some tiles of two different partitions (related by a change of coordinates, namely  by the
operator $BA$) of the hyperboloid (projected onto an open cylinder)
and the  points of the three different classes  of integer quadratic forms with $k^2-4mn=32$, lying in these tiles.

The fundamental domains are marked by  thick  black boundary:  they contain 5 points of the first orbit (circles), 5 points of the second orbit (black discs)  and 4 points  of the third orbit (squares).

Each domain in the upper regions contains 4 points of the first orbit, one point of the second, and 2 points of the third orbit. Each domain in the down regions contains 1 point of the first orbit, 4 points of the second, and 2 points of the third orbit.

The classical reduction theory introduced by Lagrange for the
indefinite forms says that there is a finite number of forms such
that $m$ and $n$ are positive and $m+n$ is less than $k$. The
reduction procedure, allowing to find these forms, can be
described in terms of the  model introduced in this work. We will
see  this relation in more detail in Part II.

The reduction theory that follows directly from our model  is in
fact closer to that expounded in the book "The sensual (quadratic)
forms" by J.H. Conway \cite{Con}, since here the `reduced' forms
are those having $mn<0$.
We prefer  this definition  for the  following  reason:  the  number  of
reduced forms by Lagrange is  equal to the number $n_u$  of forms in each domain
of the upper regions in our partition, whereas the number of reduced forms
in our  definition is  the number $n_u+n_d$ of reduced forms in the  fundamental domain.

However, we point out that the essential
new element, with respect to the known theories, is the
geometrical view-point, allowing to see the action of the group in
the space of forms, exactly as for for the modular group action on
 the Lobachevsky disc.

We introduce here also the classification of the types of
symmetries of the classes of forms.  This classification is
closely related to the classification of the symmetries of the
quadratic surds  continued fractions periods, answering more
recent questions by Arnold \cite{Ar3}, as we will show in the
second part of the article.

We will see there also how to calculate   the number  of points in
each domain for every class of hyperbolic forms  from the
coefficients of a form belonging to that class.

I am deeply grateful to Arnold who posed the problem in
\cite{Ar1}.

\section{The space  of forms and their classes}\label{sec1}

We will use also the following coordinates in the space of form
coefficients:
\[ \begin{array}{l} K  = k, \\ D  = m-n, \\ S  = m+n.
\end{array} \]

{\it Remark 1.} A point with integers coordinates $(K,D,S)$
represents a form if and only if $D=S \mod 2$.

{\it Remark 2.} The discriminant in the new coordinates reads
\[    \Delta= K^2+ D^2-S^2. \]

{\bf Definition.} A  point with integer coordinates $(m,n,k)$ or
with integer coordinates $(K,D,S)$ such that $D=S \mod 2$ is
called a {\it good point}. It will be indicated by a bold letter.

\subsection{Action of SL$(2,\mathbb{Z})$ on the form coefficients}$  \ \ \ \ \ $

Let  ${\bf f}$ be the triple $(m,n,k)$ of the coefficients of  the
form (\ref{form}), and
${\bf f}'$  the triple $(m',n',k')$, corresponding to the form
$f'$ obtained from $f$ by the action of an operator $L$ of ${\rm
SL}(2,\mathbb{Z})$. I.e., if  ${\bf v}=(x,y)$, we define $f'({\bf
v})=f(L({\bf v}))$.

We thus associate to $L$ the operator $T_L$ acting on
$\mathbb{Z}^3$ in this way:
\begin{equation}\label{ope} {\bf f'}= T_L {\bf f}.  \end{equation}

This defines an homomorphism from ${\rm SL}(2,\mathbb{Z})$ to
${\rm SL}(3,\mathbb{Z})$:  $L\mapsto T_L$.  We denote by
$\mathcal{T}$ the image of this homomorphism.  The subgroup
$\mathcal{T}$ is isomorphic to  ${\rm PSL}(2,\mathbb{Z})$, since
$T_L= T_{-L}$.

{\bf Definition.} The {\it orbit}\footnote{Words "orbit" and
"class" are synonymous.} of a good point ${\bf f}$ is the set of
points obtained applying to ${\bf f}$ all elements of the group
$\mathcal T$. The class of ${\bf f}=(m,n,k)$ is denoted by $C(\f)$ or by $C(m,n,k)$.

The following statements are obvious or easy to prove:

-- All points of an orbit are good.

-- All points of an orbit belong to the hyperboloid $K^2 + D^2
-S^2 = \Delta$. Moreover, in the elliptic case,  the orbit lies
entirely either on the upper or on the lower sheet of the
hyperboloid; in the parabolic case, it lies entirely either on the
upper or on the lower cone.

-- Every good point  belongs to one orbit.

-- Different orbits are disjoint.

\subsection{The subgroups $\mathcal{T^+}$  and  $\mathcal{T^-}$}\label{AB}$  \ \ \ \
\ $

Consider the following generators of the group ${\rm
SL}(2,\mathbb{Z})$:
\begin{equation}\label{gen}
{\bf A}=\left( \begin{array}{cc} 1 & 1 \\ 0  & 1
\end{array}\right), \ \ \ \
{\bf B}=\left(  \begin{array}{cc}  1& 0 \\ 1 & 1
 \end{array}\right), \ \ \ \
  {\bf R}=\left(  \begin{array}{cc}  0 &  1 \\ -1 & 0
  \end{array}\right).
 \end{equation}
and their inverse  operators  denoted by $\bar{\bf A}$, $\bar{\bf
B}$  and $\bar{\bf R}$.

Note that \begin{equation} \label{relAB}  {\bf R}=\bar{\bf B}{\bf
A} \bar{\bf B}={\bf A} \bar{\bf B}{\bf A} \quad {\rm and} \quad
\bar{\bf R}=\bar{\bf A}{\bf B} \bar{\bf A} ={\bf B} \bar{\bf A}
{\bf B}.\end{equation}

We denote the corresponding operators  $T_{\bf A}$, $T_{\bf B}$,
$T_{\bf R}$   of $\mathcal T$, obtained by eq. (\ref{ope}), by
$A$, $B$, $R$  and  their inverse by $\bar A$,  $\bar B$ and $\bar
R$.

{\it Remark.} The matrices of $A$ and $B$  are one the transpose
of the other\footnote{The matrices of the generators of
$\mathcal{T}$, $A$, $B$, and $R$ are, in coordinates $(m,n,k)$:
\[ A=\left(
\begin{array}{ccc}  1& 0& 0  \\ 1 & 1  & 1 \\ 2& 0& 1   \end{array} \right), \ \ \ \
B=\left( \begin{array}{ccc} 1&1 & 1
\\  0& 1  &0  \\ 0& 2& 1
\end{array} \right),\ \ \ \  R=\left( \begin{array}{ccc}  0& 1&0   \\  1& 0 &0    \\ 0&0 &
-1    \end{array} \right).  \]    The matrices of the same
generators  in coordinates $(K,D,S)$ are
\[ A=\left(
\begin{array}{ccc}  1& 1& 1  \\ -1 & 1/2  & -1/2 \\ 1& 1/2& 3/2   \end{array} \right), \ \ \ \
B=\left( \begin{array}{ccc} 1&-1 & 1
\\  1& 1/2  &1/2  \\ 1& -1/2& 3/2
\end{array} \right),  \ \ \  R=\left( \begin{array}{ccc}  -1& 0&0   \\  0& -1 &0    \\ 0&0
& 1    \end{array} \right).  \] } as well as those of  ${\bf A}$
and ${\bf B}$, whereas the transpose of ${\bf R}$ is equal to
${\bf R}^{-1}$. Since  the transpose of $R$ is equal to $\bar
R=R$, relations (\ref{relAB}) become:
\begin{equation} \label{relAB2}   R=\bar{
B}{  A} \bar{  B}={  A} \bar{  B}{  A}=\bar{ A}{ B} \bar{  A} ={
B} \bar{ A} { B}.\end{equation}

{\bf Definition.} We call $\mathcal{T}^+$ ($\mathcal{T}^-$) the
multiplicative semigroup of the elements of $\mathcal T$ generated the
identity and by the operators $A$ and $B$ ($\bar A$ and $\bar B$).

\begin{lem} \label{lem1}   a) Every operator $T\in \mathcal T^+$
($T\in \mathcal T^-$) is written in a unique way as a product of
its generators.   b) Every operator $T\in \mathcal T$ can be
written as $V S U$, where $S$ belongs to $\mathcal T^+$  and the
operators $U$ and $V$ are equal either to the identity operator or
to $R$. Statement (b) holds as well replacing $ \mathcal T^+$ by
$\mathcal T^-$.
\end{lem}

{\it Proof.}  a) There are no relations involving only operators
${\bf A}$ and ${\bf B}$ in $\SL(2,\mathbb{Z})$, hence we have no
relations involving  only $A$ and $B$. b) Relations (\ref{relAB2})
allow to transform any word in $A$, $B$, $R$ and their inverse
operators into a word of type $VSU$. This lemma is illustrated by
Figure   \ref{cap4n}; indeed,  there is an one-to-one
correspondence between the elements of the group and the domains
of the Lobachevsky disc. The element of the group corresponding to
a  domain indicates that the fundamental  domain ($I$)  is  sent
to this domain by that element. One sees that any domain in the right
half-disc is attained, from the fundamental domain, by an
operator which can be written as an element of $\mathcal T^+$
followed by the identity or by $R$, and any domain in the  left half-disc  by an
element of $\mathcal T^-$ followed by the identity or by $R$. The
multiplication by  $R$ at left acts as a reflection with respect
to the centre. Hence every domain in the right half-disc can be attained by the
operator corresponding to the domain symmetric to it with respect
to the centre, multiplied at left by $R$, and vice versa. Hence
any domain can be written using an element of $\mathcal T^+$ (as
well as  $\mathcal T^-$), multiplied at left and at right by the
identity of by  $R$.\hfill $\square$


\subsection{Symmetries of the form classes}

We   introduce in  this section some different types of symmetries that the classes of  forms may possess.

We define each one of these symmetries as the  invariance of the class under the reflection
with respect to some plane or some axis through the center of the coordinate system,  plane or axis which is {\sl not
invariant} under  the  action of the group $\mathcal T$.    So,  a priori  these symmetries
could hold no more   in another system  of coordinates.
However,  we prove  that  the  action of the group  $\mathcal T$  {\sl  preserves} each one of the symmetries,
and hence the same definitions  of the symmetries  {\sl  hold in any system of coordinates obtained by  a $\mathcal T$
coordinate transformation}. This is equivalent to say that  a  symmetry of a class of forms is  the
 symmetry  by respect to  all infinite planes (or axes), which are the images under $\mathcal T$ of
one of such symmetry planes (or axes).

\begin{figure}[h]
\centerline{\epsfbox{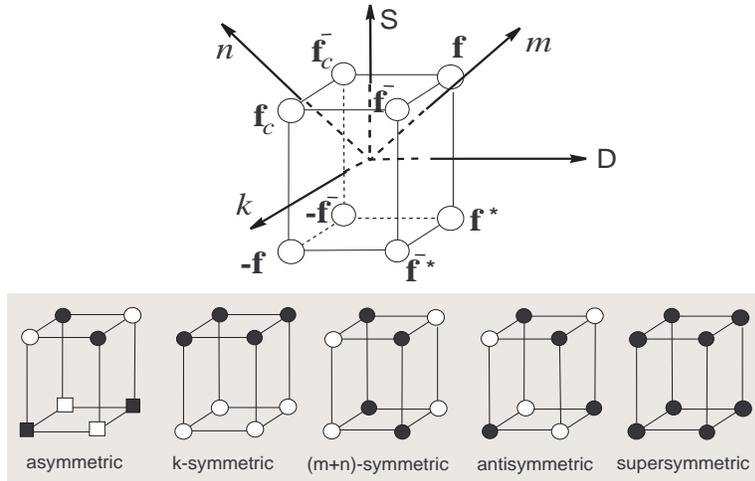}} \caption{For every symmetry type,
the forms denoted by the same symbol and the same color belong to
the same class.}\label{capp1}
\end{figure}

To every form ${\bf f}=(m,n,k)$
there correspond 8 forms, obtained by 3 involutions (see Figure
\ref{capp1}):
\begin{equation}\label{invol}
{\bf f}_c=(n,m,-k), \quad  \overline {\bf f}=(m,n,-k), \quad {\bf
f}^*=(-n,-m,k). \end{equation}

All these involutions commute, since  correspond to changes of
sign of some of the coordinates $K,D,S$. In these coordinates,
\[
{\bf f}_c=(-K,-D,S), \quad  \overline {\bf f}=(-K,D,S), \quad {\bf
f}^*=(K,D,-S).\]

The  8  forms defined by these involution on the form ${\bf f}$ lie on the same
hyperboloid as ${\bf f}$.

The {\it complementary} form ${\bf f}_c $  always belongs to the
class of ${\bf f}$,  since $\f_c=R \f$, and $R\in \mathcal T$. Note  that the
complementary form $\f_c$ of the form $\f$  satisfies  in the $(x,y)$ plane  $f_c(x,y)=f(y,-x)=f(-y,x)$
and the corresponding $\PSL(2,\Z)$-change of coordinates
is a rotation by $\pi/2$.

 The  complementary of the {\it conjugate} form $\overline \f_c=(n,n,k)$ of the form $\f=(m,n,k)$  corresponds
to the  reflection of the plane $(x,y)$  with respect to the diagonal (whose operator $(^0_1 \ ^1_0))$ does not belong
to $\SL(2,\Z)$).

 Note that the opposite form $-{\bf
f}=(-m,-n,-k)=(-K,-D,-S)$ is the complementary of the adjoint of
${\bf f}=(m,n,k)=(K,D,S)$, ($-{\bf f}={\bf f}_c^*$).

The  forms obtained from a form ${\bf f}$ by the conjugation and/or adjunction  may or may not
belong to $C({\bf f})$.   However, {\sl  if  a class contains a pair of forms related by  some involution, or a
form  which is invariant  under some involution, then the entire class is invariant under such involution}:

\begin{prop}  Let $\sigma$ be one of the following involutions:  $\sigma({\bf f})=\bar \f$ or
$\sigma({\bf f})=\f^*$ or $\sigma({\bf f})=\overline \f^*$   Suppose that, for some $\f$,  $\sigma({\bf f})\in C({\bf f})$.  Then
$\sigma({\bf g})\in C({\bf f})$ for all $\g\in C(\f)$.  \end{prop}

{\it Proof.} We have to prove the following lemma.

\begin{lem} \label{lemma-invol} The following  identities hold:
\begin{equation}\label{ident} \begin{array} {l l l}
1)   &
 {(A  {\bf f})}^* = \bar B {\bf f}^*;  &    {(B{\bf f})}^* =
\bar A {\bf f}^* ; \\  2)  &
    \overline {A  {\bf f}} = \bar A \ \overline {\bf f};
    &
 \overline{B{\bf f}} = \bar B \ \overline {\bf f}  ; \\
3)  &
    (\overline {A  {\bf f}})^* = B \overline {\bf f}^*;
    &
 (\overline{B{\bf f}})^* = A \overline {\bf f}^* . \end{array}
 \end{equation}
\end{lem}

{\it Proof of the lemma.}

Let ${\bf f}=(m,n,k)$. We have  $A{\bf
f}=(m,m+n+k,k+2m)$,    $ B{\bf f}=(m+n+k,n,k+2n)$.

1)  Since ${\bf f}^*=(-n,-m,k)$,  we have

$(A{\bf f})^*=(-m-n-k,-m,k+2m)$ and  $\bar B{\bf
 f}^*=(-n-m-k,-m,k+2m)$;

  $(B{\bf f})^*=(-n,-n-m-k,k+2n)$ and   $\bar A( {\bf f}^*)=(-n,-m-n-k,k+2n)$.

2) Since $\overline{\bf f}=(m,n,-k)$,  we have

 $\overline{A{\bf f}}=(m,m+n+k,-k-2m)$ and   $\bar A \ \overline {\bf f}=(m,m+n+k,-k-2m)$;

$\overline{B{\bf f}}=(m+n+k,n,-k-2m)$ and $\bar B \ \overline {\bf
f}=(m,m+n+k,-k-2m)$.

3) Since $\overline {\bf f}^*=(-n,-m,-k)$,  we have

 $\overline{A{\bf f}}=(m,m+n+k,-k-2m)$, $(\overline {A {\bf
f}})^*=(-n-m-k,-m,-k-2m)$ and  $B \overline {\bf
f}^*=(-n-m-k,-m,-k-2m)$;

 $\overline{B{\bf
f}}=(m+n+k,n,-k-2n)$, $(\overline {B {\bf
f}})^*=(-n,-n-m-k,-k-2n)$ and $A \overline {\bf
f}^*=(-n,-m-n-k,-k-2n)$.  \hfill $\square$

We observe now that every operator  $T\in \mathcal T$  can be written
as product  of the generators $A$, $B$ and their inverse.  Then the above lemma
implies  that  $\sigma (T {\bf f})=T'\sigma({\bf f})$ for some $T'\in \mathcal T$. Therefore, if $\sigma({\bf f})\in C({\bf f})$,
then, for every $T \in \mathcal T$,    ${\bf g}=T{\bf f}\in C({\bf f})$  and hence $\sigma{\bf g}\in C({\bf f})$.
\hfill $\square$

By  the above Proposition, the following definitions  hold in all
coordinate systems  (equivalent under $\mathcal T$-transformations).

{\bf  Definition.}

A class of forms is  said to be (see  Fig. \ref{capp1})
\begin{enumerate}

\item {\it asymmetric},  if it is only invariant under reflection with
respect to the axis of the coordinate $S$ ($K=0$, $D=0$), so
containing only pairs of complementary forms;

\item {\it $k$-symmetric}, if it is not supersymmetric but it is
invariant under reflection with respect to the plane $k=0$
($K=0$). It contains, with every ${\bf f}$,  its  {\sl conjugate} form
$\overline {\bf f}$;

\item {\it $(m+n)$-symmetric}, if it is not supersymmetric but it is
invariant under reflection with respect to  the plane $m+n=0$
($S=0$).  It contains, with every  ${\bf f}$,  its {\sl adjoint}  form
${\bf f}^*$;

\item {\it antisymmetric}, if it is not supersymmetric but it is
invariant   under reflection with respect to the planes $m=0$ and
$n=0$  ($|S|+|D|=0$). It contains, with every ${\bf f}$, its {\sl
antipodal} form (the conjugate of the adjoint) $\overline {\bf
f}^*=(-n,-m,-k)=(-K,D,-S)$.

\item {\it supersymmetric},  if it contains all 8 forms obtained by the
3 involutions;

\end{enumerate}

{\it Remarks.} 1) Note that the opposite form $-{\bf
f}$, being the   the complementary of the adjoint of the form
${\bf f}$,     belongs to the class of ${\bf f}$ only if
the class is $(m+n)$-symmetric or supersymmetric (see Fig. \ref{capp1}).

2) A  class of elliptic forms containing $\f$ cannot contain neither  $-\f$, nor $\f^*$
 nor $\overline \f^*$. Hence it is either asymmetric or  $k$-symmetric.

\section{Elliptic forms}\label{elliptic}

In this section the classification of positive definite forms is
treated in order to introduce some notions and terms which will be
used in Sections \ref{para} and \ref{hyper}.

We define a map from one sheet of the two-sheeted hyperboloid to
the open unitary disc, which gives explicitly the one-to-one
correspondence between the  integer points of an orbit on the
hyperboloid and the domains of the classical Poincar\'e model of
the Lobachevsky disc.

Let $\mathcal{P}$ be the following {\it normalized} projection
from the upper sheet of the hyperboloid $K^2+D^2-S^2=\Delta$
($\Delta<0$) to the disc of unit radius. Let ${\bf p}=(K,D,S)$ be
a point on the hyperboloid (see Figure  \ref{capp2h}, left), ${\bf
p}'$  its projection to the disc of radius $\rho=\sqrt{-\Delta}$
from the point $O'$,  and $\mathcal{P}{\bf p}$ the image of the
normalized projection. We have
\begin{equation}\label{tilde} \mathcal{P}{\bf p}= \Big\{ \begin{array} {c} \widetilde
K=\frac {K}{\rho+S} \\ \widetilde D= \frac{D}{\rho+S}. \end{array}
\end{equation}

Let $L=\begin{pmatrix} a & b \\ c& d
\end{pmatrix} $ be any operator of ${\rm SL}(2,\mathbb{Z})$ and
$T_L$ its corresponding operator defined by  (\ref{ope}).

{\bf Definition.} We define the operator $\widetilde L$ acting on
the disc of radius $1$ by: \begin{equation}\label{elle} \widetilde
L ( \mathcal{P}{\bf p})=  \mathcal{P} (T_L {\bf p}).
\end{equation}

On the other hand, the operator $L\in  {\rm SL}(2,\mathbb{Z})$
defines another map  from the disc to itself.  Indeed,  let $H_L$
be the homographic operator acting on the upper complex half-plane
$\{z\in\mathbb{C}: \Im(z)\ge 0 \}$:
\begin{equation}\label{hom} H_L z = \frac{az+b}{cz+d}.
\end{equation}

\begin{figure}[h]
\centerline{\epsfbox{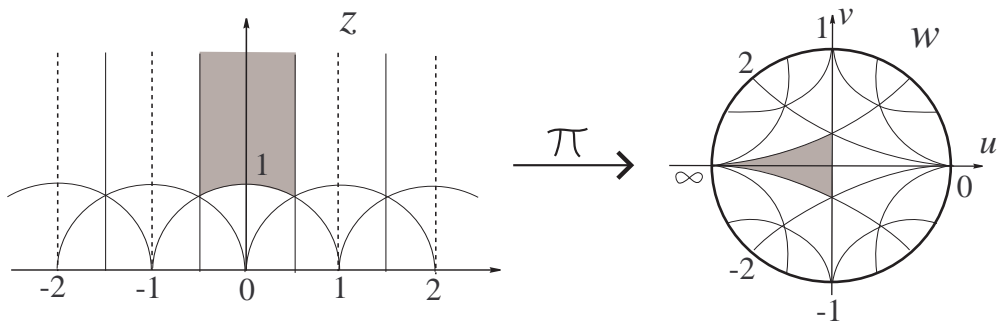}} \caption{The fundamental domain is
gray}\label{capp3}
\end{figure}

The following map $\pi: z \rightarrow w$   sends the upper complex
half plane to the unitary complex disc $\{w\in\mathbb{C}: |z| \le
1 \}$ (Figure \ref{capp3})
\begin{equation} \label{pi}  w=  \pi z = \frac{1+iz}{1-iz}.
\end{equation}

  Define the operator $\widehat L$ acting on the
complex unitary disc by:
\begin{equation}\label{hat} \widehat L( \pi z)=  \pi( H_L z ).
   \end{equation}

\begin{prop} \label{pro1} The actions of operators $\widehat L$  and $\widetilde L$ on
the unitary disc do coincide under the identification:
  \[    \widetilde D = \Re w,  \ \ \  \   \widetilde K =  \Im w. \]
\end{prop}

The proof is done by writing  explicitly the operators
corresponding to the ${\rm SL}(2,\mathbb{Z})$ generators and by
comparing their actions in coordinates. \hfill $\square$

{\it Remark.} The group of operators defined by eq. (\ref{elle}) is
isomorphic to the homographic group of the operators $H_L$, and to
the group $\mathcal {T}$, i.e., to ${\rm PSL}(2,\mathbb{Z})$: we
will denote its generators by the same letters as the
corresponding generators of $\mathcal T$.

\begin{figure}[h]
\centerline{\epsfbox{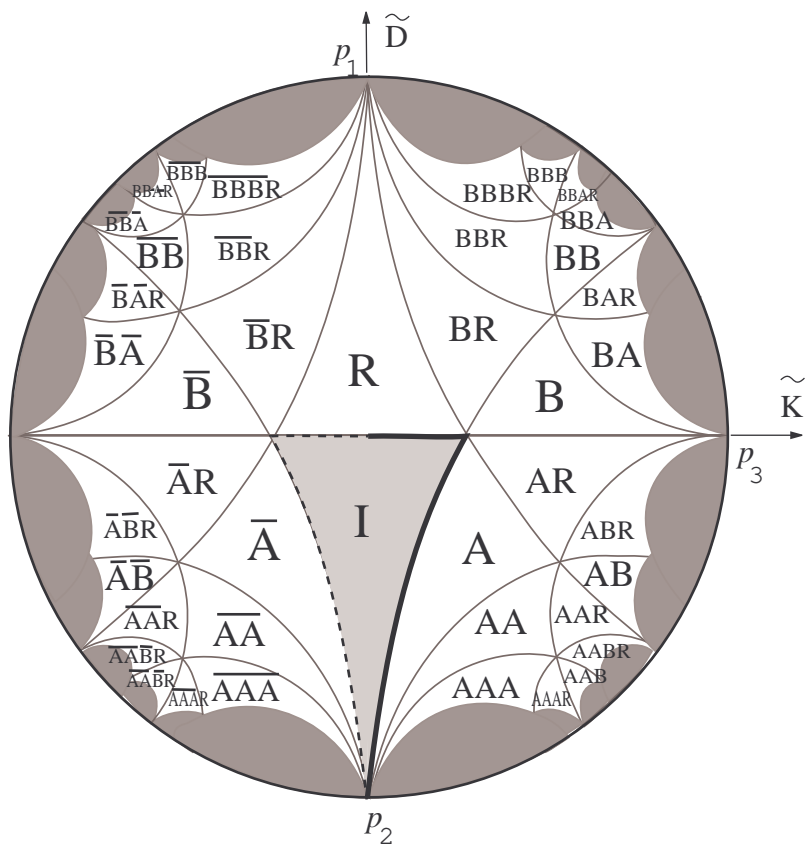}} \caption{ A finite set of domains
in the Lobachevsky disc with coordinates $\widetilde K, \widetilde
D$.  } \label{cap4n}
\end{figure}

We take as coordinates in the Lobachevsky disc the pair
$(\widetilde K,\widetilde D)$. Hence our Lobachevsky disc is
obtained from the unitary complex disc  with coordinate $w=u+iv$
(see Figures \ref{capp3} and \ref{cap4n})  by the reflection with
respect to the diagonal $v=u$.


In Figure \ref{cap4n}   the Lobachevsky disc  is shown with a
finite sets of domains. The fundamental domain  is indicated by
the letter I (the bold line at the frontier belong to it, the
dotted line does not).   The other domains are attained applying
to I the elements of $\PSL(2,\Z)$, written in terms of $R$,
$A,B$ and their inverse.

Because of relations (\ref{relAB2}) involving these generators,
the expressions are not unique. We have chosen this representation
in order to see the meaning of Lemma \ref{lem1}, which will be
determinant in Section \ref{hyper}.

{\it Remark.}  The choice of the fundamental domain is arbitrary.
Every domain  has the shape of a triangle with one and only one
corner on the circle at infinity.  Choose  an arbitrary triangle,
labeled by the operator $M$, and orient it as the fundamental one,
with the corner on the circle at bottom. Then the triangle
adjacent at right is labeled by $MA$,  the triangle adjacent at
left by $M\bar A$, and the triangle adjacent at the base by $MR$,
and so on. Replacing $M$ by $I$ in all triangles,  the arbitrarily
chosen triangle becomes the new fundamental domain.

\begin{figure}[h]
\centerline{\epsfbox{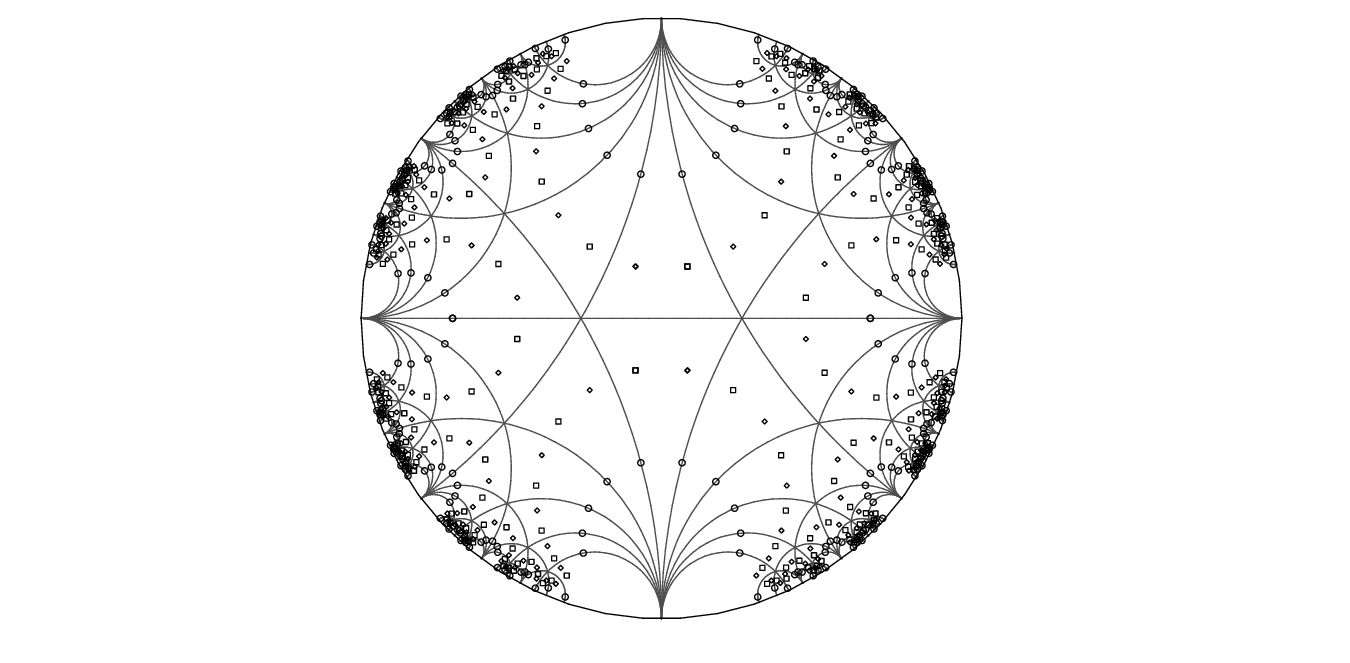}} \caption{Finite subsets of the 3
distinct orbits in the case $\Delta=-31$ (the representative point
($m,n,k$) is in the fundamental domain). Two asymmetric orbits:
boxes (2,4,-1) and rhombi (2,4,1),  and one $k$-symmetric: circles
(1,8,1). Unitary disc coordinates are $(\widetilde K, \widetilde D
)$. } \label{capp4}
\end{figure}

Every orbit has one and only one point in every domain. In Figure
\ref{capp4} the Lobachevsky disc with a finite subset of domains
is shown together with a finite part of the three distinct orbits
in the case $\Delta=-31$.

{\it Remark.} The opposite, the antipodal, and the adjoint  of a
positive (negative) definite quadratic form  $f$ are negative
(resp., positive) definite quadratic forms, so they cannot belong
to the same class of $f$. For this reason, a class of elliptic
form may have only two types of symmetries: it is either
$k$-symmetric or asymmetric.


\subsection{The hierarchy of the points at the
infinity}\label{hier}$  \ \ \ \ \ $

This section is important  for the study of hyperbolic forms.

The points of the circle $C$ (the circle at infinite bounding the
Lobachevsky disc) having rational coordinates,  inherit, by the
${\rm PSL}(2,\mathbb{Z})$ action, a hierarchy,  that we will
recall. These points will be thus distinguished into {\it points
of zero, first, second (and so on) generation}. We notice that
this hierarchy   will be the starting point of our
 partition of
the de Sitter world.

Denote by $\pi'$ the  composition of  mapping $\pi$ (see eq.
(\ref{pi})) with the reflection of the disc-image ($|w|\le 1$)
with respect to the diagonal ($\Im w =\Re w$).

We consider, because of the symmetry of the picture, only the
right  part of the circle $C$. This semicircle ($\widetilde
K\ge0$) is the image by $\pi'$ of the half real line $x\equiv\Re
z>0$ of the half-plane where the homographic operators act.

The points $p_1=(\widetilde K=0,\widetilde D=1)$ and
$p_2=(\widetilde K=0,\widetilde D=-1)$ of zero generation (i.e.,
the end-points of this semicircle) are the images by $\pi'$ of the
points $x_1=0$ and $x_2=\infty$ of the real half-line.

The rational points $x_i$ of the half real line are written as
fractions, i.e.: $0\equiv0/1$, $\infty\equiv 1/0$, $q\equiv q/1$,
if $q\in \Z$ etc,  and points $p_i$ are their images by $\pi'$ on
the circle $C$.

 Here  A and B are  the generators of the homographic
group that correspond  to the generators ${\bf A}$ and ${\bf B}$
of SL$(2,\mathbb{Z})$. \begin{equation}  {\rm A}: x \rightarrow
\frac{x+1}{1}
; \ \ \ {\rm B}: x \rightarrow  \frac{x}{x+1}.
\end{equation}

Consider the iterate action of such generators on the points $x_1=0$ and
$x_2=\infty$.  We have firstly: \begin{equation} \label{x1x3} {\rm A}
x_1=x_3; \ \ \ {\rm B} x_1= x_1;  \ \ \   {\rm A} x_2= x_2; \ \ \
{\rm B} x_2=x_3, \ \ \  \end{equation} where $x_3=1/1$ is the
preimage by $\pi'$ of the point $p_3=(\widetilde K=1, \widetilde
D=0)$.

{\bf Definition.} {\it  The points $x_i$ of the $n$-th generation
are obtained from the point of first generation $x_2=1/1$ applying
to it  all the $2^n$ words in the generators {\rm A} and {\rm B}.
The hierarchy and the order of these points is shown in the
following scheme, where {\rm T} indicates any word of $2^{(n-1)}$
generators: }

\centerline{\epsfbox{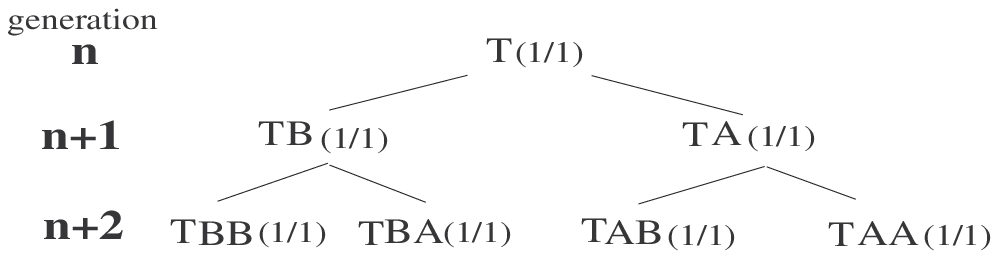}}

The points of all generations have a nice algebraic property that
we recall.

{\bf Definition.} We call {\it sons} of the point T$(1/1)$,
belonging to the  $n$-th generation, the points TA$(1/1)$ and
TB$(1/1)$ of the $(n+1)$-th generation. T$(1/1)$ is thus the {\it
father} of his sons. In the scheme above the segments indicate the
relations father-son.

{\bf Farey rule.} {\it  The coordinate of a point of the $n$-th
generation can be calculated directly from those of his father and
of the nearest (i.e. the closest in the line where these points
lie) ancestor to his father, by the rule shown in the following
scheme:}

\centerline{\epsfbox{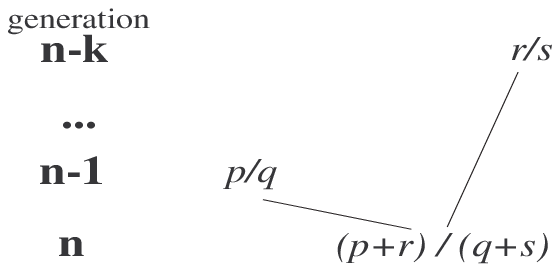}}

In the following scheme  of the hierarchy every point $x_i$ is
connected by a segments to its  2 sons, to its father and to its
closest ancestor. Note that the descendants of B$(1/1)$ are the
inverse fractions of the descendants of A$(1/1)$.

\centerline{\epsfbox{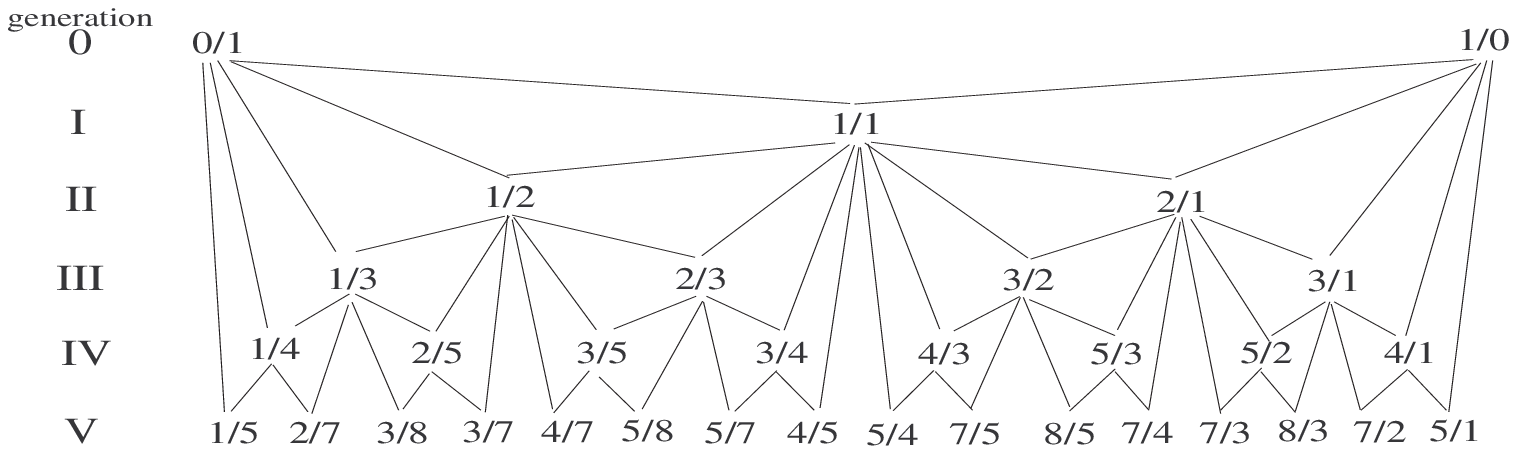}}

{\it Remark. } By this procedure  all positive rational numbers
are covered.

The relations of order and the hierarchy are preserved by the map
$\pi'$. Hence we have the same ordering and hierarchy on the
rational points  on the circle $C$.

Note that the point $(\widetilde K=-1,\widetilde D=0)$ on the
circle $C$ is the image by $\pi'$ of point $-1/1$: the above
construction can be repeated by the iterate action of the inverse
generators A$^{-1}$ and B$^{-1}$ on point $(-1/1)$, by considering
the point $-1/0=-\infty$ as  the preimage of  point $p_2$. In this
way one finds  the rational points in the half real line $x<0$,
and the ordering and the hierarchy of the points with rational
coordinates of the left part of  circle $C$.

The point $x_i=\frac{p}{q}$ is sent  by $\pi'$ to the point of $C$
with coordinates $(\widetilde K, \widetilde D)$:
  \begin{equation} \label{pith} \widetilde K= \frac
  {2pq}{p^2+q^2},\ \ \  \ \ \widetilde D=
  \frac{p^2-q^2}{p^2+q^2}.\end{equation}

{\it Remark.} The map $\pi'$    defines a map from the set of
rational numbers to the set of Pythagorean triples $\{(a,b,c)\in
\mathbb{Z}^3: \ a^2+b^2=c^2 \}$:
\[  \frac{p}{q} \rightarrow  \left( a=2pq, \ \ \ b=q^2-p^2, \ \ \  c=p^2+q^2 \right).  \]

\section{Parabolic forms}\label{para}

{\bf Definition.} A good point $(K,D,S)$ is called {\it
Pythagorean} if $K^2+D^2=S^2$.  If $K,D,S$ have no common
divisors, then the triple is said {\it simple}. If $(K,D,S)$ is
Pythagorean, then the set of points $\{ (\lambda K, \lambda D,
\lambda S), \ \lambda \in \mathbb{Z} \}$, all Pythagorean, is
called {\it Pythagorean line}.

Pythagorean points belong to the cone $\Delta=0$, and any good
point belonging to the cone is Pythagorean.

\begin{lem} \label{lem3} There exists an one-to-one correspondence between the
Pythagorean lines and  the points with rational coordinates $p_i$
on the circle $C$. \end{lem}

{\it Proof.} By the last remark of the preceding section,  to
every point $p_i$ on the circle $C$  we  associate a Pythagorean
triple. This triple represents a good point because $b-c=0 \mod
2$. On the other hand, having a simple good point $(K,D,S)$, the
equations
\[    K=2pq, \ \ \ D=p^2-q^2, \ \ \  S=p^2+q^2  \]
have solution $p=\sqrt{(S+D)/2}$; $q=\sqrt{(S-D)/2}$, i.e.,
$p=\sqrt m$, $q=\sqrt n$.

Since $S=m+n$ and $D=m-n$, and in this case $K=2\sqrt{mn}$, $m$
and $n$ have no common divisors, otherwise triple $(K,D,S)$ should
be not simple. But the equality $K^2=4mn$, with $m$ and $n$
relatively prime, implies that $m=p^2$ and $n=q^2$ for some
integers $p$ and $q$. Hence to every simple Pythagorean triple we
associate a point $p_i$ on the circle $C$ and vice versa. The
Pythagorean line corresponding to $p_i$ is the line through 0 and
$p_i$ on the cone of parabolic forms.  \hfill $\square$

\begin{thm} \label{the1}  There are infinitely many classes of forms with $\Delta=0$.
In particular,  all forms of type $ax^2$, ($m=a,n=0,k=0$) with
$a\in \Z$, belong to different orbits, and every orbit on the cone
contains a form of this type.
\end{thm}

{\it Proof.} We prove that for all $a\in \mathbb{Z}$ the orbits
containing the points $(K=0,D=a,S=a)$ are distinct.  Suppose that
point ${\bf r}=(0,b,b)$ belong to the same orbit of point ${\bf
p}=(0,a,a)$. Thus the form $f=ax^2$  is in the same class of the
form $f'=bx^2$. This means that there exists an operator
$L=\begin{pmatrix} \alpha & \beta \\ \gamma  & \delta
\end{pmatrix}$ of $\SL(2,\Z)$ such that $a(\alpha x+\beta y)^2=bx^2$. This
can be satisfied only by $\beta=0$, and, since
$\alpha\delta-\gamma\beta=1$, by $\alpha=\delta=1$.  Hence $a=b$.
 On the other hand, any parabolic form $f=mx^2+ny^2+kxy$, being
$k^2-4mn=0$, can be written as $a(\alpha x+\beta y)^2$, $a$ being
the greatest common divisor of $(m,n,k)$. For every pair of
integers $\alpha$ and $\beta$, there exist two integers $\gamma$
and $\delta$ such that $\alpha\delta-\gamma\beta=1$.  So the
inverse of the operator $L$ is the operator of $\SL(2,Z)$ transforming
form $f$ into $a x^2$.
 \hfill $\square$



\section{Hyperbolic forms}\label{hyper}

\begin{figure}[ht]
\centerline{\epsfbox{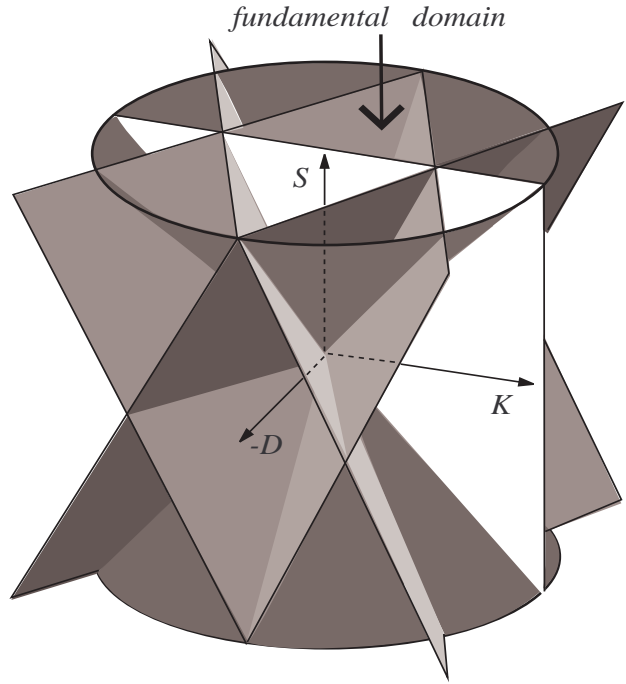}} \caption{The  fundamental domain
  in the space of forms coefficients. }\label{capp19}
\end{figure}

Let $X$ be the infinite set of planes  through the origin in the
3-dimensional space ($K,D,S$), obtained from the plane  $D=0$ by
the action of group $\mathcal T$.  These planes  subdivide the
interior of the cone ($K^2+D^2 < S^2$) into domains. (Some of
these planes are shown in Figure \ref{capp19}). These planes,
intersecting both sheets of the two-sheeted hyperboloid,
subdivides them into domains;  those which belong to the upper
sheet, are are projected by $\mathcal{P}$ (see Section
\ref{elliptic}) to the interior of domains of the Lobachevsky
disc.

The closure $\overline X$ of $X$ contains also the planes tangent
to the cone along all Pythagorean lines.

The intersection of the planes of $\overline X$ with the
one-sheeted hyperboloid  $H$ ($K^2+D^2-S^2=1$), at the exterior of
the cone, forms a net of lines which is dense in $H$. For this
reason there is no a natural prolongation of the Poincar\'{e}
model of the Lobachevsky disc to the exterior of the disc - i.e.,
on the de Sitter world.

The horizontal section $S=1$ of the set $\overline X$ gives, at
the interior of the unit circle, the Klein model of the
Lobachevsky disc. The arcs of circles between two points of the
circle at infinite of the Poincar\'e model are substituted by the
chords connecting these points. We are interested in the
prolongations of these chords outside the disc: the description of
the de Sitter world is based on a particular subset of these
chords: they are the "limit" chords (the section $S=1$ of the
planes tangent to the cone) i.e., the tangents to the circle at
all rational points of it, provided with the hierarchy explained
in Section \ref{hier}.

\subsection{The Poincar\'e model of the de Sitter world}$  \ \ \ \ \ $

In analogy with the standard  projection $\mathcal P$ from the
upper  sheet of the two sheeted  hyperboloid to the Lobachevsky
disc, I have chosen the following mapping $\mathcal Q$ from the
hyperboloid $H$, with equation $K^2+D^2-S^2=\Delta$ in coordinates
$(K,D,S)$, to the open cylinder $C_H$:
\[ C_H=\{ (K,D,S) \ :\  K^2+D^2=1,  \ \ |S|<1 \}. \]

\begin{figure}[h]
\centerline{\epsfbox{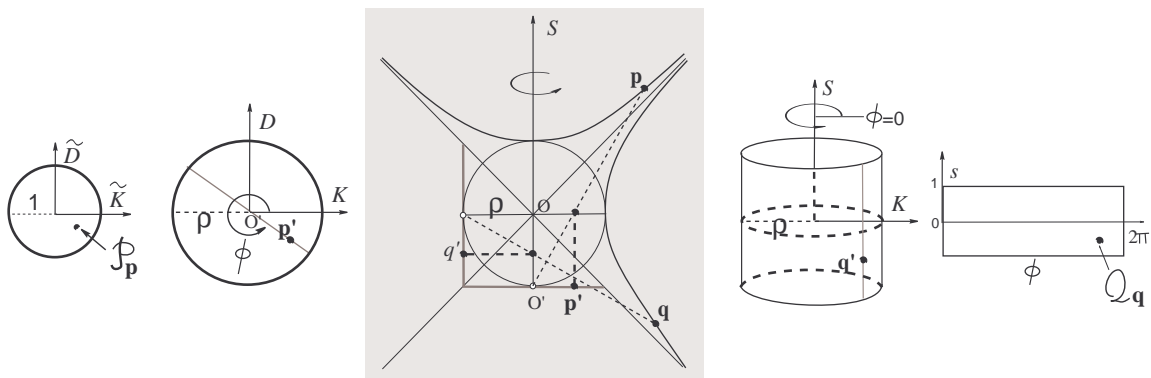}} \caption{Projection $\mathcal{P}$,
left, and $\mathcal{Q}$, right  }\label{capp2h}
\end{figure}

{\bf Definition.}  Let $\rho=\sqrt{\Delta}$ and
$r=\sqrt{K^2+D^2}$. The coordinates on the cylinder $C_H$ are:
\begin{equation}\label{Q} s=\frac{S}{r+\rho},   \end{equation} and $\phi$ is the angle  defined by
the relations: $K=r \cos \phi$ and $D=r\sin \phi$ (see Figure
\ref{capp2h}, right).

The border of the cylinder consists of two circles, denoted $c_1$
($s=1$)  and $c_2$ ($s=-1$).

\begin{figure}[h]
\centerline{\epsfbox{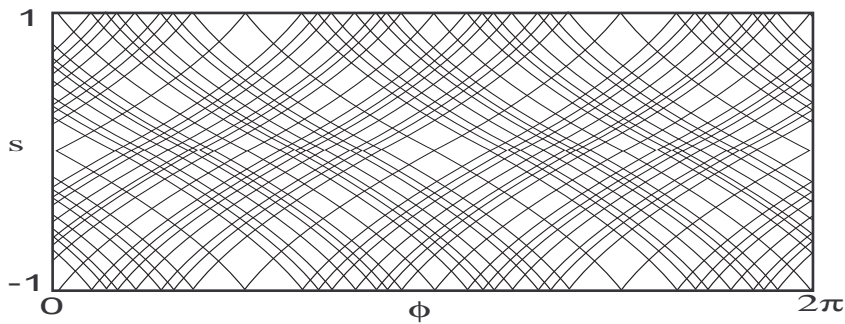}} \caption{Projection on the
cylinder $C_H$ of few lines of the infinite set of lines,
intersection of the planes tangent to the cone along Pythagorean
lines. }\label{capp5}
\end{figure}

{\it Remark.}  The lines of intersection of the planes tangent to
the cone with the hyperboloid $H$ are straight lines, being
generatrices of the hyperboloid.

{\bf Definition.} We denote by $H^0$ and  by $H_R^0 $ the domains
\[  H^0= \{(K,D,S) \in H : |S|<|D| ,  D>0  \}; \]
\[  H_R^0= \{(K,D,S) \in H : |S|<|D| ,  D<0  \}. \]

{\it Remark.} $H^0_R=R H^0$.

Since the mapping $\mathcal Q$ is one-to-one, for simplicity we
denote by the same letters the domains on $H$ and their images
under $\mathcal Q$ on the cylinder $C_H$.

Circles $c_1$ and $c_2$ on the planes $S=1$ and $S=-1$ coincide
with the circle $C$ at the infinite of the Lobachevsky disc. Hence
also on $c_1$ and $c_2$ the points having rational coordinates
($\widetilde K, \widetilde D$) are mapped, according to
(\ref{pith}), into the points $p_i$ of first, second, third  -and so
on- generation with the hierarchy explained in Section \ref{hier}.

\begin{figure}[h]
\centerline{\epsfbox{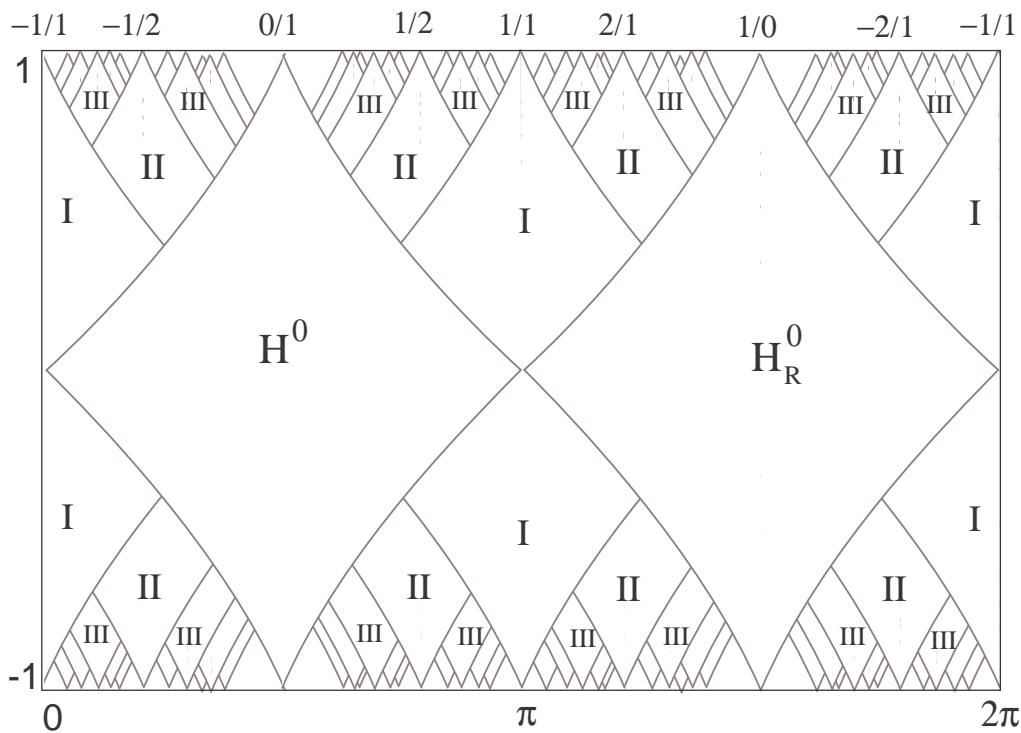}} \caption{The Poincar\'e model of
the de Sitter world: domains of  first, second and third
generations are indicated by I,II,III}\label{capp6}
\end{figure}

In Figure \ref{capp6}, on the upper circle $c_1$, the points $p_i$
up to the second generations are denoted   by the corresponding
rational numbers $x_i$.

Note that the upper (and lower) vertices of the domain $H^0$ and
$H^0_R$, with coordinates  $\phi=\pi/2$ and $\phi=3\pi/2$,
correspond to the points $x_1=0/1$ and $x_2=1/0=\infty$.

{\bf Definition.} We denote by $H^{x_i}$ and $H^{-1/x_i}$ the
domains obtained by $H^0$ and $H^0_R$ by a rigid translation to
right, such that the the upper   vertex of $H^0$ transfers to
points $x_i$, and the upper vertex of $H^0_R$ transfers to point
$-1/x_i$.  So, domains $H^{x_i}$ inherit the hierarchy of points
$x_i$: $H^0=H^{0/1}$ and $H^0_R=H^{1/0}$ are called {\it rhombi of
zero-generation}, $H^{1/1}$ and $H^{-1/1}$ {\it rhombi of first
generation},  $H^{1/2}$, $H^{-2/1}$, $H^{-1/2}$, $H^{2/1}$ {\it
rhombi of second generation} and so on.

{\bf Definition.}  We denote by $H^O$,  $H^I$, $H^{II}$, $H^{III}$
etc the unions of the rhombi of zero, first, second, third etc
generations.

{\bf Definition.} We call the {\it Poincar\'e model} of the
hyperboloid $H$ the cylinder $C_H$ provided with the subdivision
into domains obtained by the following procedure. Let
$\mathcal{H}^0=H^0 \cup H^0_R$ be the {\it domain of zero
generation}. Let $\mathcal{H}^I=H^I \setminus (H^I \cap \overline
{H^O})$  the domain of first generation, $\mathcal{H}^{II}= H^{II}
\setminus ( H^{II}\cap (\overline {H^O} \cup \overline{H^I} ))$
the domain of second generation and so on (where $\overline{H^{\bf
n }}$ is the closure of $H^{\bf n }$): the domain of the $n$-th
generation are thus obtained as
\[   \mathcal{H}^{\bf n}=H^{\bf n}
\setminus (H^{\bf n}\cap (\overline {H^O} \cup \overline{H^I}
\cup\overline{H^{II}}\cup \cdots \cup \overline{H^{\bf n-1}})).
\]

Figure  \ref{capp6} shows the domains of different generations.
Note that the domain of $n$-th generation, $\mathcal{H}^{\bf n}$,
has $2^{n+1}$ connected components.

\begin{figure}[h]
\centerline{\epsfbox{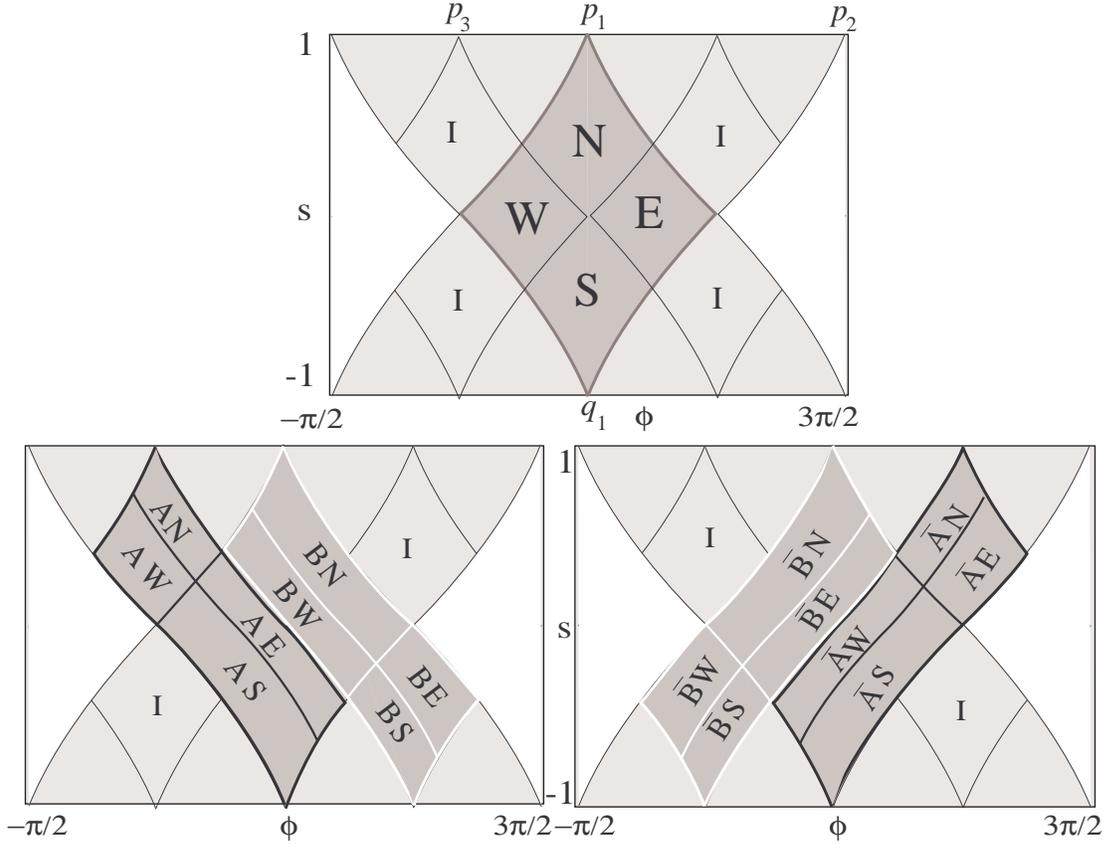}} \caption{Images by $A$, $B$, $\bar
A$ and $\bar B$ of $H^0$. Letter ${\rm I}$ indicates the four
connected components of the domain of first generation $\mathcal
H^I$}\label{capp7}
\end{figure}

The action of operators $A,B$ and their inverse on the domain
$H^0$ is shown in  Figure \ref{capp7}, where $\phi$ varies from
$-\pi/2$ to $3\pi/2$, so that domain $H^0$ is at the center of the
figure.

We have subdivided $H^0$ into sub-domains, denoted by N,S,E,W.
The operators  $A, B$ and their inverses  map $H^0$ partially to
itself, and partially outside $H^0$.

{\it Remark.} The corresponding actions on $H^0_R$ are obtained
taking into account the following relations among the generators,
coming from (\ref{relAB2}):
\begin{equation} \label{relaz} AR= R \bar B;\ \ \  BR= R \bar A; \ \ \
\bar A R= R B; \ \ \ \bar B R= R A.   \end{equation}

\begin{lem} \label{lem4} The parts of images of $H^0$ (of  $H^0_R$) by $A,B$ and
their inverse which are not in $H^0$  (in $H^0_R$) are disjoint
and cover the first generation domain $\mathcal H^I$. \end{lem}

{\it Proof.} It suffices to calculate the action of the generators
on the vertices of $H^0$, and on the intersections of the frontier
lines of it with the frontiers of the rhombi of   first
generation. Note that the action of $A$ and $B$ on the points at
infinite (on the circles $c_1$ and $c_2$) is given by equations
(\ref{x1x3}), remembering that any  point $q_i$  on the circle
$c_2$, symmetrical  of the point $p_i$ on $c_1$, as opposite
vertex of the same rhombus, is the opposite ($q_i=-p'_i$) of point
$p'_i$ on the circle $c_1$ at distance $\pi$  from $p_i$. So, for
instance, the image under $A$ of the extreme north $p_1$ of $H_0$
is $p_3$, according to (\ref{x1x3}), whereas the image under $A$
of the extreme south, $q_1$, is  $q_1$, since $q_1=-p_2$, and
$Ap_2=p_2$ (see Figure \ref{capp7}). \hfill $\square$

\begin{figure}[h]
\centerline{\epsfbox{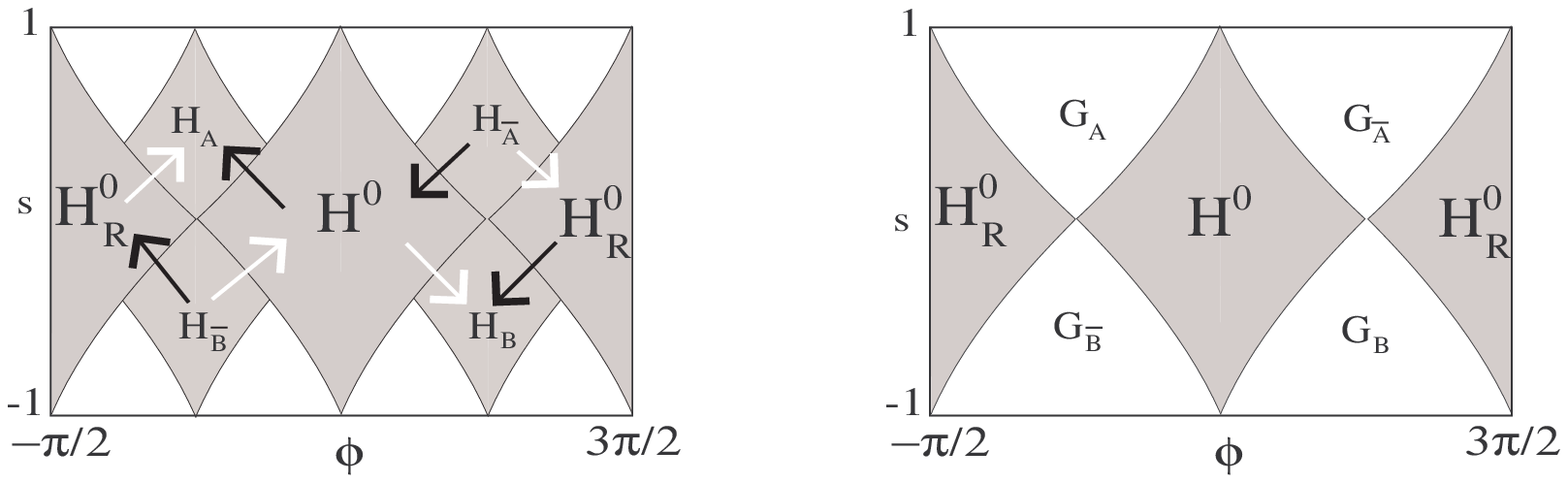}}\caption{}\label{capp8}
\end{figure}

{\bf Definition.} We denote the four connected components of
$\mathcal H^I$, where $|S|<|K|$ and $|S|>|D|$,  by (see Figure
\ref{capp8},left):
\begin{equation}\label{HAB}
\begin{array}{ll}
           H_A & =  \{ (K,D,S)\in \mathcal H^I: S > 0, K >0 \};  \\
 H_{\bar A} & =  \{ (K,D,S)\in  \mathcal H^I: S > 0, K < 0 \}; \\
 H_B & =  \{ (K,D,S)\in  \mathcal H^I: S< 0, K <0 \}; \\
  H_{\bar B} & =  \{ (K,D,S)\in \mathcal H^I: S < 0, K >0 \}.
 \end{array} \end{equation}

{\it Remark.} Note that the subscript  index of a connected
component of the first generation domain $\mathcal H^I$ indicates
the operator mapping one part of $H^0$ into this connected
component. Moreover
\[       H_{\bar A} = R H_A;
  \ \ \ \  H_{\bar B}= R H_B. \]

Figure \ref{capp8}, left,  shows the four connected components of
$\mathcal{H}^I$ reached from  $H^0$ and from $H^0_R$ by means of
the operators $A$ (black arrow) and $B$ (white arrow). The reversed
arrows must be read as the inverse operators.

{\bf Definition.} The domain $G^0 \equiv H\setminus
(\overline{H^0\cup H^0_R})$ consists of four disjoint connected
components (see Figure \ref{capp8}, right), denoted by:
\[ \begin{array}{ll}
           G_A & =  \{ (K,D,S): S > |D|, K >0 \};  \\
 G_{\bar A} & =  \{ (K,D,S): S > |D|, K < 0 \}; \\
  G_B & =  \{ (K,D,S): S < -|D|, K <0 \}; \\
   G_{\bar B} & =  \{ (K,D,S): S < -|D|, K >0 \}.
 \end{array} \]

{\it Remark.}  $G_{\bar A}=RG_A$, $G_{\bar B}=R G_B$.

\begin{thm} \label{the2} There is an one-to-one correspondence between the
operators of  $\mathcal T^+$ and the connected components of the
domains of all generations of order $r>0$ inside $G_A$  (inside
$G_B$). There is an one-to-one correspondence between the
operators of $\mathcal T^-$ and the connected components of the
domains of all generations of order $r>0$ inside $G_{\bar A}$
(inside $G_{\bar B}$).

In fact, for all $T\in \mathcal T^+$ that are  the product of $n$
generators (of type $A$ and $B$), the domains $TH_A$, $TH_B$,
$\bar{T}H_{\bar A}$, $\bar{T}H_{\bar B}$ cover all the connected
components of the domain of the $n$-th generation,
$\mathcal{H}^{\bf n+1}$, belonging respectively to $G_A$, $G_B$,
$G_{\bar A}$, $G_{\bar B}$.
\end{thm}

\begin{figure}[h]
\centerline{\epsfbox{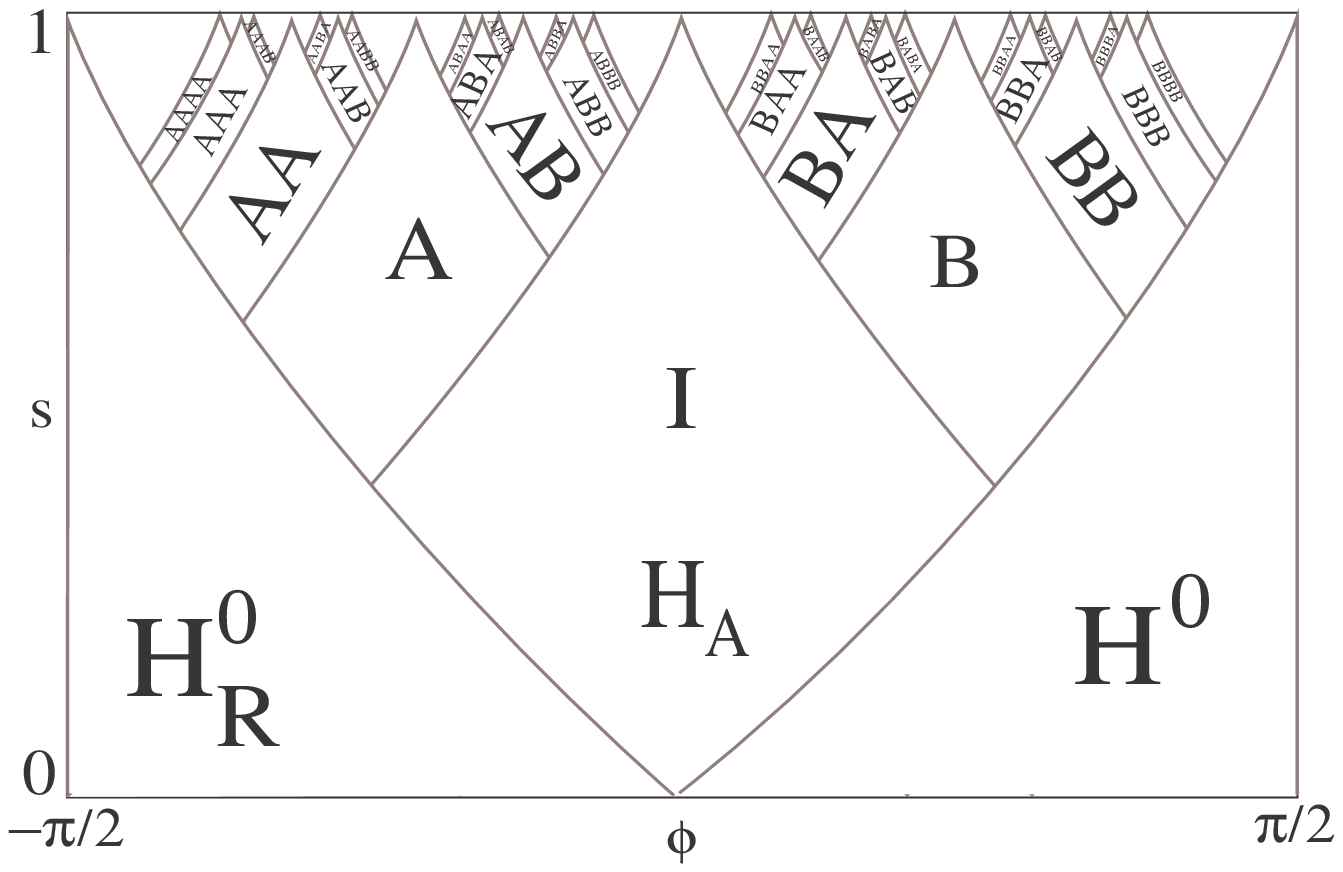}}\caption{}\label{capp9}
\end{figure}

The proof of this theorem is  a computation as well. The
correspondence between the operators of $\mathcal T^+$ and the
domains in $G$   is shown in Figure \ref{capp9}. \hfill $\square$

{\it Remark.} Theorem \ref{the2} implies that {\it domains $H_A$
and $H_B$ behave as fundamental domains for the action of
$\mathcal T^+$} in $G_A$ and $G_B$ respectively, whereas {\it
domains $H_{\bar A}$ and $H_{\bar B}$ behave as fundamental
domains for the action of $\mathcal T^-$ in $G_{\bar A}$ and
$G_{\bar B}$ respectively.}


\subsection{Coordinates changes}$  \ \ \ \ \ $

Here we see how  the Poincar\'e model of the de Sitter world changes under a
change of coordinates in the space  of the coefficients corresponding to a
$\SL(2,\Z)$ changes of coordinates in the plane  where  the forms are defined.

As we have said  in the introduction, the situation in the de Sitter world is different from
that of the Poincar\'e model  in the Lobachevsky  plane, where a change of coordinates
by an operator $L \in \SL(2,\Z)$ in the plane
corresponds to replace the fundamental  domain by its image by  $T_L$, which is another domain of
of  the tiling.

The procedure in the de Sitter world  is  complicated  by the fact that  the images  under $\mathcal T$
 of the domains overlap each other,  and, as we have seen in the preceding section, the key of
 the model  is to choose,  for all non-fundamental domains, the  parts  of them  which do not overlap
 under  the action of  the semigroups $\mathcal T^+$ or $\mathcal T^-$.

Let us consider  an operator $T\in \mathcal T$,  operating our change of coordinates. The  image by
$T$ of $H^0$ and $H^0_R$  are   the fundamental domains of the new tiling. In the following figure  $T=\bar A$.

\centerline{\epsfbox{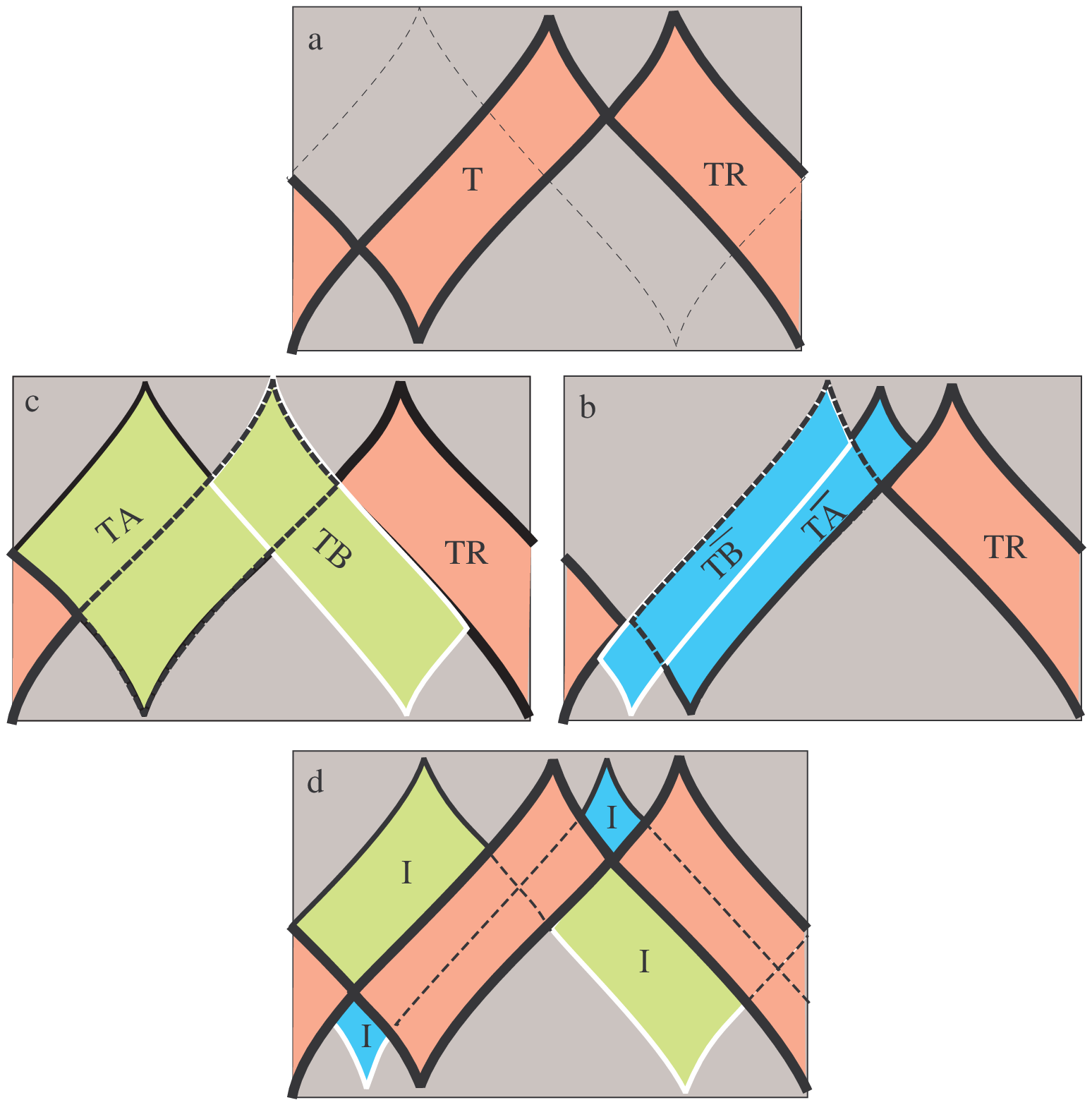}}

In part  (a), the domains marked by $T$ and $TR$   are  the image by  $T$ of $H^0$ and $H^0_R=RH^0$.
They represents  the new {\it rhombi of zero generation}. (The  dotted lines  show the boundaries of $H^0$ and $H^0_R$).

Part (b)  of the figure shows   the images by  $T\bar A$ and by $T\bar B$ of the fundamental domain  $H^0$  (marked by  $T\bar A$ and $T\bar B$ respectively) and  part (c)  of the figure shows   the images  by  $TA$ and by $TB$ of the fundamental domain  $H^0$ (marked by  $T A$ and $T B$ respectively).
The union of these  images form the  {\it rhombi of the first generation}.

Part (d). To obtain the {\it domain of first generation} we have to exclude from these rhombi the parts of them which overlap with the rhombi of   zero  generation (which  become  the fundamental domains $H^0$ and $H^0_R$).  In this way we obtain 4 disjoint components (denoted by I), which are the image by $T$ of the domains $H_A$, $H_{\bar A}$, $H_B$ and $H_{\bar B}$ of the standard model introduced in the previous section.
The  procedure to  build the tiling continues  analogously  to that already explained.


\subsection{Hyperbolic orbits}$  \ \ \ \ \ $

\begin{figure}[h]
\centerline{\epsfbox{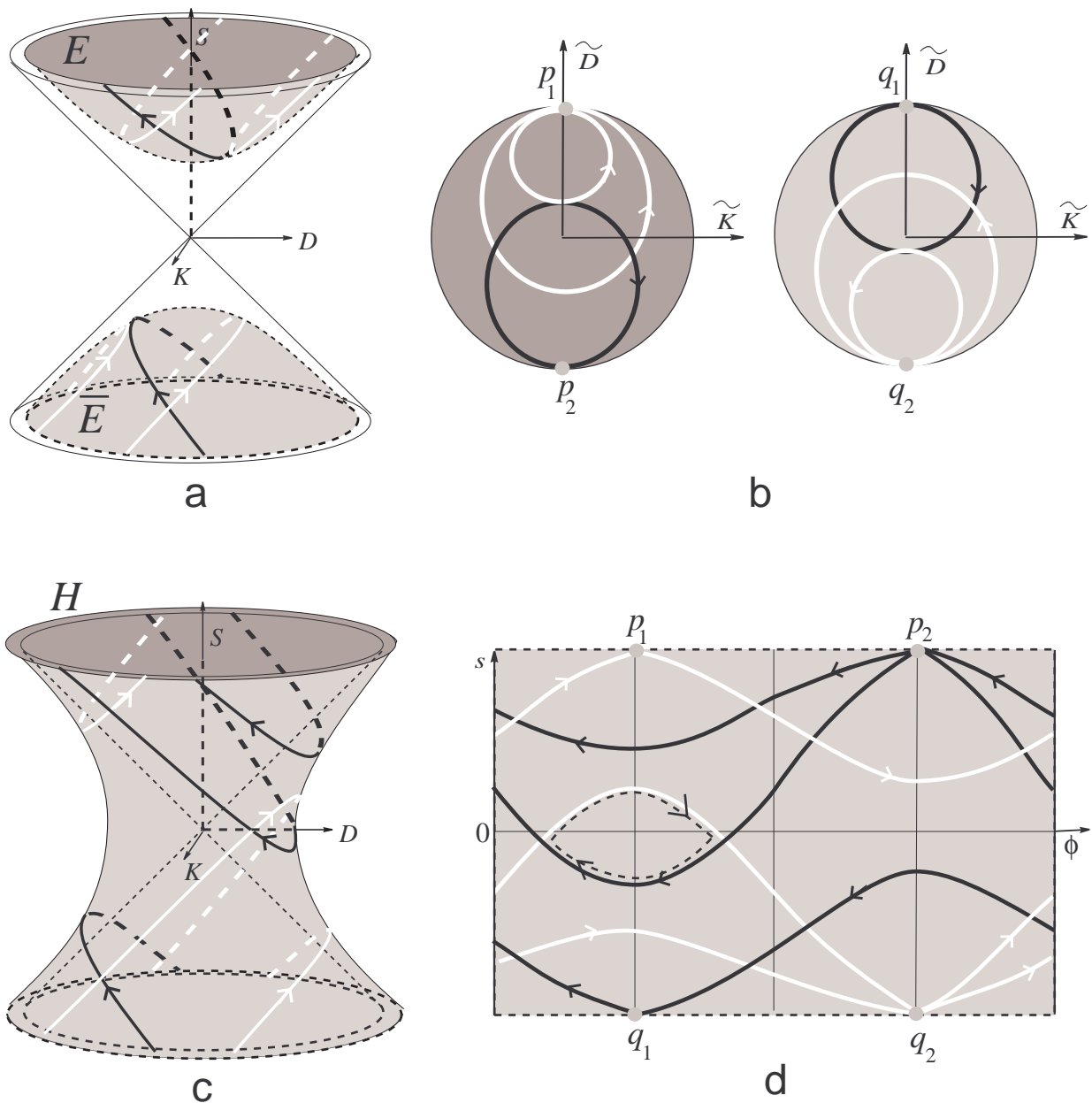}} \caption{a) Neck-laces on the
two-sheeted hyperboloid $\{E,\bar E\}$;  b) their projection on
the Lobachevsky disc; c)  Neck-laces on the one-sheeted
hyperboloid $H$; d) their projection on the cylinder
$C_H$}\label{capp10}
\end{figure}

Figure \ref{capp10} illustrates the reason of the difference
between the action of $\mathcal T$ on the two-sheeted hyperboloid
$\{E,\bar E\}$ and on the one-sheeted hyperboloid $H$, as we will
explain.

Consider a generic good point ${\bf f}$ on a hyperboloid  in the
space $(K,D,S)$.

{\bf Definition.} The sequences of points  $[A^j {\bf f}]$ and
$[B^j {\bf f}]$, where $j\in \mathbb{Z}$, are called   {\it
neck-laces} and are denoted by $\omega_{\bf f}(A)$ and
$\omega_{\bf f}(B)$, respectively. The ordering of $\mathbb{Z}$
provides any neck-lace with an orientation.

\begin{prop} \label{pro2} Every neck-lace $\omega_{\bf f}(A)$ (every
neck-lace $\omega_{\bf f}(B)$) lies on the intersection of the
hyperboloid with a plane of the family $S= - D + \alpha$
(respectively, $S= D +\alpha$), $\alpha\in \mathbb{Z}$. \end{prop}

{\it Proof.}  For every ${\bf f}=(K,D,S)$, the vector $(A -I){\bf
f}$ is orthogonal to the vector $(0,1,1)$ (the normal vector to
the family of planes $S+D-\alpha=0$) and the vector $(B-I){\bf f}$
is orthogonal to the vector $(1,-1,0)$ (the normal vector to the
family of planes $S-D-\alpha=0$). \hfill $\square$

{\it Remark.} The notion of neck-lace is independent of the type
(elliptic or hyperbolic) of orbits. Figure \ref{capp10} shows some
lines where the neck-laces lie on  both two-sheeted and
one-sheeted hyperboloids.  Black colour is used for neck-laces  of
type $\omega_{\bf f}(A)$ and white colour for neck-laces of type
$\omega_{\bf f}(B)$. The arrow shows the orientation of the
neck-lace.

Note that in the elliptic case  the projections of the neck-laces
lie on   horocycles tangent to the circle $C$ (the $\infty$) at
the points $\widetilde K=0$ ($p_1$ and $p_2$ in Figure
\ref{capp10}{\bf b}).

{\bf Definition.} We call {\it semi-orbits} $O^+_{\bf f}$ and
$O^-_{\bf f}$ the sets of good points  reached  by the action of
the semigroups $\mathcal{T}^+$ and $\mathcal{T}^-$, respectively,
excluding the identity,  on a good point ${\bf f}$.

For every ${\bf f}$ in $E$,  $O^+_{\bf f}$ never contains ${\bf
f}$. This can be seen starting by any point  in $E$, and  trying
to reach it by a path composed of  pieces of neck-laces always in
the same direction of the arrow  (i.e., by an operator either of
$\mathcal{T}^+$, or of $\mathcal{T}^-$). In $H$ the situation is
different: a series of good points obtained one form the previous
one applying in sequence either operator $A$ or operator $B$
(i.e., forming a path composed of pieces of neck-laces in the
positive direction) can lie on a cycle (see for instance the
dotted-line in Figure \ref{capp10}{\bf d}).

We will consider separately the case  when $\Delta$ is a square
number.

\subsection{The case of $\Delta$ different from a square number.}$  \ \ \ \ \ $

\begin{thm} \label{the3}  For every $\Delta$ different from a square number there
exists some good point ${\bf f}$  on the hyperboloid $H_\Delta$
and an operators $T$ of $\mathcal T^+$ ($\mathcal T^-$), such that
$T{\bf f}={\bf f}$. Such a point belongs to $H^0$.\end{thm}

{\it Proof.} The proof  follows from  some
lemmas.

\begin{lem} \label{lem5} The domain $H^0$ ($H^0_R$) contains  a finite number of good
points. \end{lem}

{\it Proof.} Since $mn=(S^2-D^2)/4$, the domain $H^0$  contains  all
forms where $m>0$ and $n<0$, whereas $H^0_R$ contains all forms
where $m<0$ and $n>0$.  From the definition of $\Delta$ we obtain
\[   4mn= k^2  -\Delta.  \]
Since in $H^0$ the product $mn$ is negative, the above equality is
fulfilled by a finite number of values of $k$, less than
$\sqrt{\Delta}$. For each one of these values a finite number of
product $|4mn|$ equals $\Delta-k^2$.   \hfill $\square$

Note that the condition $\sqrt{\Delta} \notin \mathbb{Z}$ means
that $m$ and $n$ cannot vanish, hence we  have $|D|\not=|S|$. Thus
$G^0$ ($G^0\equiv H \setminus (\overline{H^0 \cup H^0_R})$)
contains all hyperbolic forms where $m$ and $n$ have the same
sign.

{\bf Definition.} A {\it cycle of length $t$} ($t>1$) is a sequence of
points $[{\bf f}_1,{\bf f}_2,\dots, {\bf f}_t]$  such that ${\bf
f}_i=T_{i-1}{\bf f}_{i-1}$ ($i=2,\dots,t$) and ${\bf f}_1=T_t {\bf
f}_t$, where each one of operators $T_1, T_2, \dots, T_t$ is
either $A$ or $B$. A cycle of length $t$ is indicated by
$\gamma_{\bf f}(T_1,\dots,T_t)$ where ${\bf f}={\bf f}_1$. An
equivalent notation of the cycle $\gamma_{\bf f}(T_1,\dots,T_t)$
is, evidently, $\gamma_{\bf g}(T_j,\dots,T_t,T_1,\dots T_{j-1})$,
where ${\bf g}={\bf f}_j$.

\begin{lem} \label{lem6} Every  good point ${\bf f}$ in $H^0$  satisfies:
$A {\bf f} \in  H^0$  iff  $B {\bf f}  \in  H_B$;   $B {\bf f} \in
H^0$ iff $A {\bf f}  \in  H_A$, and, moreover,  $\bar{A} {\bf f}
\in H^0$ iff $\bar{B} {\bf f} \in H_{\bar B}$;  $\bar{B} {\bf f}
\in H^0$ iff $\bar{A} {\bf f} \in H_{\bar A}$. Analogous
statements hold substituting $H^0$ with $H^0_R$.
\end{lem}
{\it Proof.} By the action of $A$, the point ${\bf f}$, in
coordinates $(m,n,k)$, is sent to $A {\bf f}=(m,m+n+k,2m+k)$, and,
by $B$, to $B {\bf f}=(m+n+k,n,2n+k)$. Since ${\bf f}$ belongs to
$H^0$, $m>0$ and $n<0$. Now, if $m+n+k>0$ then $A {\bf f}$ belongs
to $G^0$, and $B {\bf f}$ belongs to $H^0$, whereas if $m+n+k<0$
then $A {\bf f}$ belongs to $H^0$, and $B {\bf f}$ to $G^0$.
Similar inequalities hold for the inverse generators.  By Lemma
\ref{lem4}, the image by $A$, $B$, $\bar{A}$ and $\bar{B}$ of a
point in $H^0$, if it is in $G^0$, belongs respectively to $H_A$,
$H_B$, $H_{\bar A}$, and $H_{\bar B}$. Because of relations
(\ref{relaz}), analogous arguments hold for $H_R^0$. \hfill
$\square$

\begin{lem} \label{lem7} For every good point ${\bf f}$ of $H^0$ ($H^0_R$) there is an integer
$t>1$ such that ${\bf f}$  belongs to a unique cycle $\gamma_{\bf
f}(T_1,\dots,T_t)$. \end{lem}

{\it Proof.} Let ${\bf f}_1={\bf f}$.  By  Lemma \ref{lem6},
either $A {\bf f}_1$ or $B {\bf f}_1$ belongs to $H^0$.  Let ${\bf
f}_2=T_1{\bf f}_1$ be in $H^0$, being $T_1=A$ or $T_1=B$. Now, let
${\bf f}_3$ be the point, among $A{\bf f}_2$ and $B {\bf f}_2$,
which belongs to $H^0$, etc. We find in this way  a sequence of
points ${\bf f}_1,{\bf f}_2, {\bf f}_3,\dots$ in $H^0$. Since, by
Lemma \ref{lem5}, $H^0$ contains only a finite number of good
points, the sequence of point $[{\bf f}_1,{\bf f}_2,\dots,{\bf
f}_t]$ must be periodic, i.e., we find eventually a cycle
$\gamma_{\bf f}(T_1,\dots,T_t)$. For $H^0_R$ the proof is
analogous. \hfill $\square$

\begin{lem} \label{lem8} Different cycles are disjoint.\end{lem}

{\it Proof.} By Lemma \ref{lem6}, any point of a cycle determines
the others, hence if two cycles have a common point, then they
coincide.  \hfill $\square$

\begin{lem} \label{lem9} $H^0$ ($H^0_R$) contains at least one good point if
$\sqrt{\Delta} \notin \mathbb{Z}$.  \end{lem}

{\it Proof.}  The discriminant, $\Delta$, either is divisible by
$4$ or $\Delta = 4d+1$.  If $\Delta=4d$, then the point
$(m=d,n=-1,k=0)$ belongs to  $H^0$  and the point $(m=-1,n=d,k=0)$
belongs to  $H^0_R$. If $\Delta=4d+1$, then the point
$(m=d,n=-1,k=1)$ belongs to $H^0$ and the point $(m=-1,n=d,k=0)$
belongs to  $H^0_R$.  \hfill $\square$

We have thus completed the proof of the Theorem\footnote{An
alternative  proof of this theorem should consist in
proving that the set of eigenvectors corresponding to the
eigenvalue $\lambda=1$  of the operators of $\mathcal{T}^+$
contains an integer good vector  $(v_1,v_2,v_3)$ such that
$v_1^2+v_2^2-v_3^2=4d+e$,  for every $d\in \mathbb{N}$ and $e\in
\{0,1\} $ whenever $4d+e$ is different from a square number.}
\ref{the3}. Indeed, let ${\bf f}$ be a good point of $H^0$,
defined by Lemma \ref{lem9}. By Lemma \ref{lem7}, it belongs to a
cycle $\gamma_{\bf f}(T_1,\dots, T_t)$. Hence $T=T_t T_{t-1}
\cdots T_2 T_1 \in \mathcal T^+$ satisfies $T {\bf f}={\bf f}$.
\hfill $\square$

\begin{thm} \label{the5} There is an one-to-one  correspondence between the orbits
under $\mathcal{T}$ of the good points of $H_\Delta$ (when
$\Delta$  is different from a square number) and the cycles in
$H^0$.\end{thm}

{\it Proof.} The point $R {\bf f}$, obtained from ${\bf f}$ by a
shift by $\pi$ on the cylinder $C_H$, belongs to the orbit of ${\bf
f}$. Hence any orbit on $C_H$ is invariant under a shift by $\pi$.
All the following statements on $H^0$ hold analogously for
$H_R^0$.

By Lemma \ref{lem8}, different cycles are disjoint. Every cycle
belongs to some orbit, by definition of orbit.  We have to prove:
a) that every orbit contains a  cycle; b) that this cycle is
unique.

a) Let ${\bf f}$ be a point, non belonging to a cycle, so ${\bf
f}\in G^0$. Let us suppose that ${\bf f}\in G_A$ (in the other
cases the proof is analogous). By Theorem \ref{the2}, there exists
a unique operator $T$ of $\mathcal T^+$ such that ${\bf
g}=T^{-1}{\bf f}$ belongs to $H_A$.  Hence ${\bf h}=\bar A {\bf
g}=\bar A  T^{-1}{\bf f} $ is inside $H^0$. But if a point belongs
to $H^0$, then it belongs to a cycle by Lemma \ref{lem7}, and
hence the orbit of ${\bf f}$ contains a cycle.

b)  We have to prove that  the point ${\bf h}$ in $H^0$, obtained
from ${\bf f}$ by the above procedure, is unique, i.e., that  we
cannot reach another point of $H^0$ non belonging to the cycle
$\gamma_{\bf h}$ by an operator of $\mathcal T$.  By Lemma
\ref{lem1}, every operator of $\mathcal T$ can be written as
$USV$, where $S$ belongs to $\mathcal T^+$ or to $\mathcal T^-$,
and $V$ and $U$ are equal to the identity or to the operator $R$.
Hence we try to reach $H^0$ from ${\bf f}$  by means of operators
of such types. We have to start by $R$, reaching a point ${\bf p}$
of $G_{\bar A}$. As before, there is only one operator of
$\mathcal {T}^-$ such that  $p$ is the image by it of a point  in
$H_{\bar A}$, indeed: \[ \begin{array} {l}
  {\bf p}=R{\bf f}  \ \ \ {\bf p }\in G_{\bar A},  \quad {\rm hence }\\
  {\bf p}=R\ T {\bf g}=  \hat T  R{\bf g}  \ \ ,
  \end{array} \]
where  $\hat T $  is the operator obtained by $T$ replacing each
$A$ by $\bar{B}$ and each $B$ by $\bar{A}$, and vice versa. The
operator $\hat T$ belongs to  $\mathcal T^-$,  and,   being ${\bf
g}\in H_A$, point ${\bf j}=R{\bf g}$ is in $H_{\bar A}$. Now we
can reach   a point either in $H^0$ (by $A$), or in $H^0_R$ (by
$B$).
 Since $g=A{\bf h}$, we obtain
\[ B {\bf j} = R R \ B {\bf j} = R \ \bar{A}R \ R {\bf g}= R \ \bar A
A{\bf h}= R {\bf h}. \] Therefore in this case the point in
$H^0_R$, reached from ${\bf p}$ in $G_{\bar A}$,  is exactly $R
{\bf h}$.
 On the other hand, $A {\bf j}\in H^0$ is equal to:
\[ \begin{array}{l}
A {\bf j} = R R \ A {\bf j} = R \ \bar{B}R \  R {\bf g} = R \
\bar{B}\  A {\bf h}= A R \ A   {\bf h}= B {\bf h}. \\
  \end{array}\]
Since $A {\bf h}$ is in $H_A$, $B {\bf h}$ is in $H^0$ (by Lemma
\ref{lem6}) and belongs to the orbit of ${\bf h}$. The proof is
completed. \hfill $\square$

\begin{figure}[h]
\centerline{\epsfbox{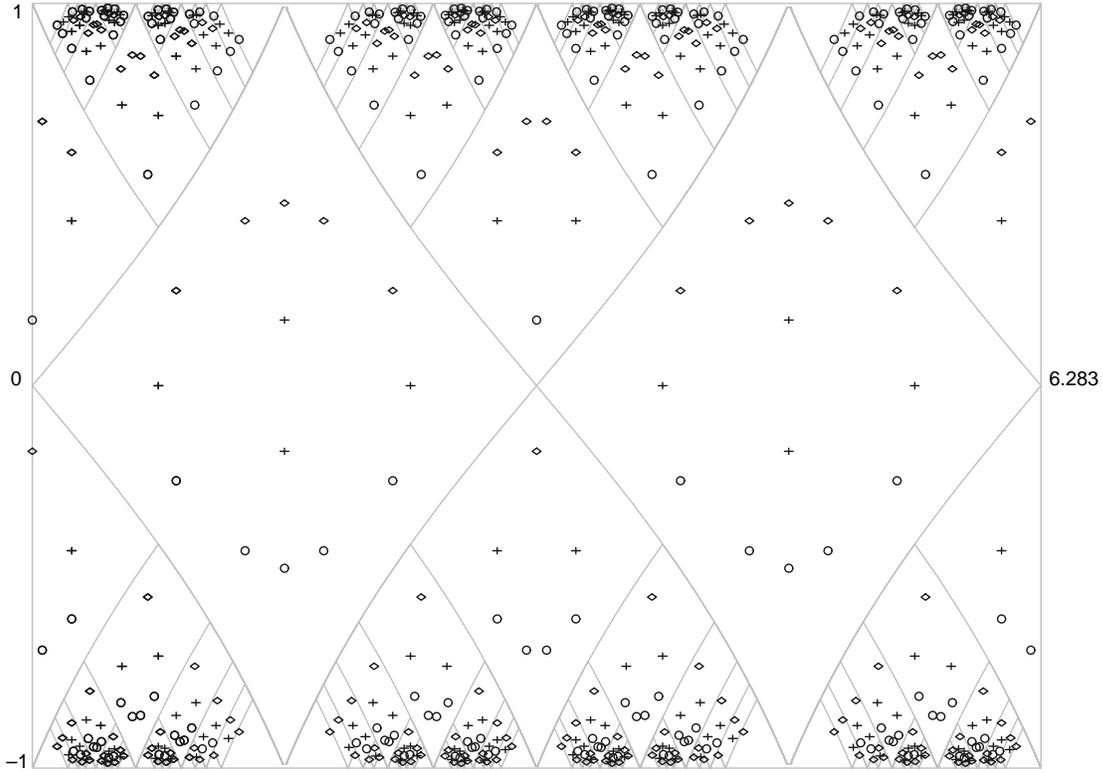}} \caption{The initial parts of the
3 distinct orbits with $\Delta=32$ projected to the cylinder. The
initial points of orbits in coordinates $(K,D,S)$  are: (-6,0,2),
circles, $k$-symmetric orbit; (6,0,2), diamonds, $k$-symmetric
orbit; (-4,4,0), crosses, supersymmetric orbit. }\label{capp12}
\end{figure}

Figure \ref{capp12}  shows  examples of orbits projected to the
cylinder $C_H$.

Theorems  \ref{the3} and \ref{the5} implies the following
two Corollaries.

\begin{cor} \label{cor1} The set of goods points in $H^0$ ($H^0_R$) is partitioned into
cycles. \end{cor}

\begin{cor} \label{cor2}
Let operator $T$ of $\mathcal T^+$ ($\mathcal T^-$) satisfy $T{\bf
f}={\bf f}$, for some good point ${\bf f}\in H^0$, and $T$ be
composed by $t$ generators of type $A$ or $B$. Then in $H_0$ the
orbit of ${\bf f}$ has  $t$ points ${\bf f}_i$, including ${\bf
f}$,  satisfying $\tilde T_i {\bf f}_i={\bf f}_i$, $\tilde T_i$
($i=1,\dots, t$) being obtained from $T$ by a cyclic permutation
of the sequence of the $t$ generators defining it. Such points
belong to $H^0$, and no other points of the same orbit belong to
$H^0$.
\end{cor}

\begin{thm} \label{the4} The operator $T$ defining the cycle in $H^0$
is the product of $t$ operators, $t_A$ of type $A$ and $t_B$ of
type $B$ ($t_A+t_B=t)$ iff:

-- in $H_A$ and  in $H_{\bar A}$, as well as in every domain in $G_A$
and $G_{\bar A}$, there are $t_B$ points;

-- in $H_B$ and in $H_{\bar B}$, as well as  in every domain in
$G_B$ and $G_{\bar B}$, there are $t_A$ points.
\end{thm}

\begin{figure}[h]
\centerline{\epsfbox{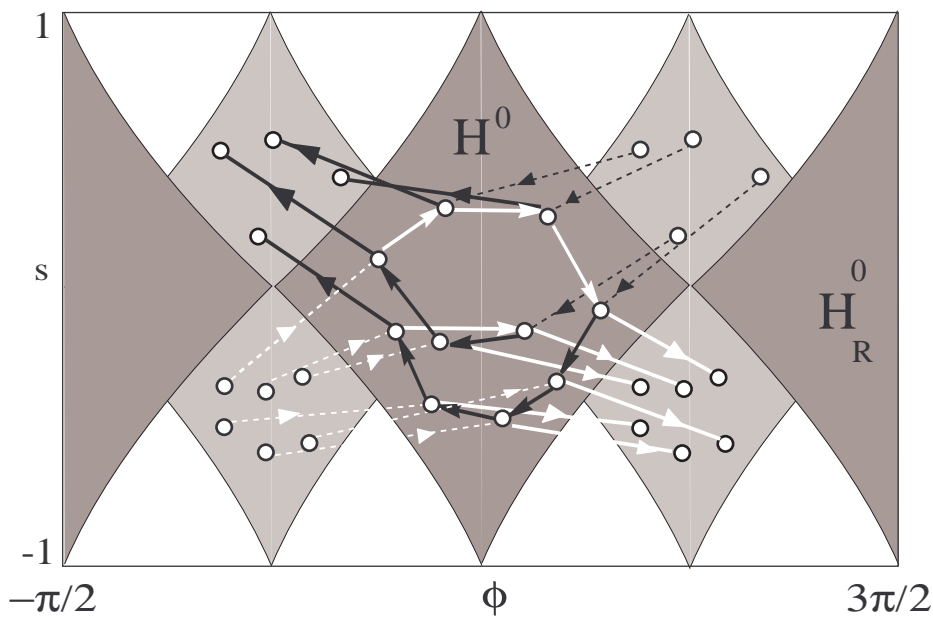}} \caption{The cycle in $H^0$ of an
asymmetric orbit, for $\Delta=624$, containing 10 points, among
which four are mapped by $A$ to $H_A$, six  by $B$ to $H_B$, four
by $\bar A$ to $H_{\bar A}$ and six by $\bar B$ to $H_{\bar B}$
}\label{cap11n}
\end{figure}

{\it Proof.} By Theorem \ref{the3}, the set of points of every
cycle in $H^0$ is subdivided into two disjoint subsets: the set of
points whose image under $A$ belongs to $H_A$   and the set of
points whose image under $B$  belongs to $H_B$. By Lemma
\ref{lem6}, the image by $A$ of a point is inside $H_A$ if and
only if its image under $B$ is inside $H^0$.  Moreover, there are
$t_B$ of such points, if and only if $T$ contains $t_B$ generators
of type $B$ (see Figure \ref{cap11n} and Figure \ref{capp8}).
Similarly, the image by $B$ of a point is inside $H_B$ if and only
if the image under $A$ is inside $H^0$; and there are $t_A$ of
such points, if and only if $T$ contains $t_A$ generators $A$.
Similarly, in virtue of the same Lemma \ref{lem6}, we obtain an
analogous statement, considering the inverse of the operator $T$:
the image by $\bar{A}$ of a point is inside $H_{\bar A}$ if and
only if the image under $\bar{B}$ is inside $H^0$; and there are
$t_B$ of such points, if and only if $\bar{T}$ contains $t_B$
generators of type $\bar{B}$, etc. (see Figure \ref{cap11n}).
\hfill $\square$

\subsection{The case of $\Delta$ equal to a square number}$  \ \ \ \ \ $

This case is essentially different from the preceding one because
the coefficients $m$ and $n$  of the form can vanish. This means
that for such forms either $S=D$ or $S=-D$, and therefore the
corresponding good points lie on the boundary of $H^0$ and of
$H^0_R$.

\begin{thm} \label{the6} On the hyperboloid $H_\Delta$, $\Delta=\rho^2$, $\rho\in
\mathbb{Z}$, there are exactly $\rho$ orbits, corresponding to the
points   $(K=\rho,\  D=r,\  S=r)$, $r=0\dots \rho-1$. \end{thm}

{\it Proof.} The proof follows from the following lemmas.

{\bf Definition.} We use the following notations for the open
segments, belonging to the common frontier of $H^0$ and  $H_A$,
$H_{\bar A}$, $H_B$ and $H_{\bar B}$:
\[
 \begin{array} {ll}
 F_A= &\{ (K=\rho,D=r, S=r), 0<r< \rho \} \ ; \\
 F_{\bar A}= &\{ (K=-\rho,D=r, S=r), 0<r< \rho \}  \ ; \\
 F_{\bar B}= &\{(K=\rho,D=r, S=-r), 0<r<\rho  \}   \ ; \\
 F_B= &\{ (K=-\rho,D=r, S=-r), 0<r< \rho \} . \\
\end{array}  \]

\begin{lem} \label{lem10} There is an one-to-one  correspondence between the images of the
set $F_A$ by the operators of $\mathcal T^+$ which are product of
$n$ generators of type $A$ and $B$  and the lower-right sides of
the frontiers of all domains of the $(n+1)$-th generation inside
$G$. There is an one-to-one correspondence between the images of
set $F_B$ by the operators of $\mathcal T^+$ (product of $n$
generators) and the upper-left sides of the frontiers of all
domains of the $(n+1)$-th generation inside $G_B$. (Similar
statements hold for $\mathcal T^-$, $F_{\bar A}$ and $F_{\bar
B}$).\end{lem}

{\it Proof.} $F_A$ is the lower-right side of the frontier of
$H_A$.  The action of $\mathcal T^+$ on $F_A$ is deduced from that
on $H_A$ (see Theorem \ref{the2}).  Note that $F_B$ is the
upper-left side of the frontier of $H_B$, etc. \hfill $\square$

\begin{lem} \label{lem11} Every  good point ${\bf f}$ in $H^0$ satisfies:

$A {\bf f} \in  H^0$  iff  $B {\bf f}  \in  H_B$;  $B {\bf f} \in
 H^0$ iff $A {\bf f} \in F_A$;

$\bar A {\bf f} \in  H^0$  iff  $\bar B {\bf f}  \in  H_{\bar B}$;
$\bar B {\bf f} \in
 H_0$ if $\bar A {\bf f} \in H_{\bar A}$.

 Moreover,    ${A} {\bf f} \in F_{ A}$ iff
${B} {\bf f} \in F_{ B}$ and  $\bar{A} {\bf f} \in F_{\bar A}$ iff
$\bar{B} {\bf f} \in F_{\bar B}$.
\end{lem}

{\it Proof.} This lemma is the version of Lemma \ref{lem6} when
$\Delta$ is equal to  a square number. Indeed, if ${\bf
f}=(m,n,k)$ and ${\bf g}=A{\bf f}\in F_A$, then $(m+n+k=0)$. This
implies that ${\bf g}'=B{\bf f}=(m+n+k,n,k+2n)$ belongs to $F_B$.
The cases of the inverse generators are similar. \hfill $\square$

 \begin{lem} \label{lem12} For every $\Delta=\rho^2$, the orbit of the point $(K=\rho,
D=0, S=0)$ is supersymmetric: it contains, with point
$(-\rho,0,0)$, all the lower points of the frontier of all domains
in $G$ and $G_{\bar A}$  and all upper points of the frontiers of
all domains in $G_B$ and in $G_{\bar B}$. These points are the
only points of the orbit. \end{lem}

{\it Proof.} This lemma is a consequence of Lemma \ref{lem4}  and
\ref{the2}, being points $(K=\pm\rho, D=0, S=0)$ the lower points
of domains $H_A$ and $H_{\bar A}$ and the upper points of domains
$H_B$ and $H_{\bar B}$. \hfill $\square$

\begin{figure}[h]
\centerline{\epsfbox{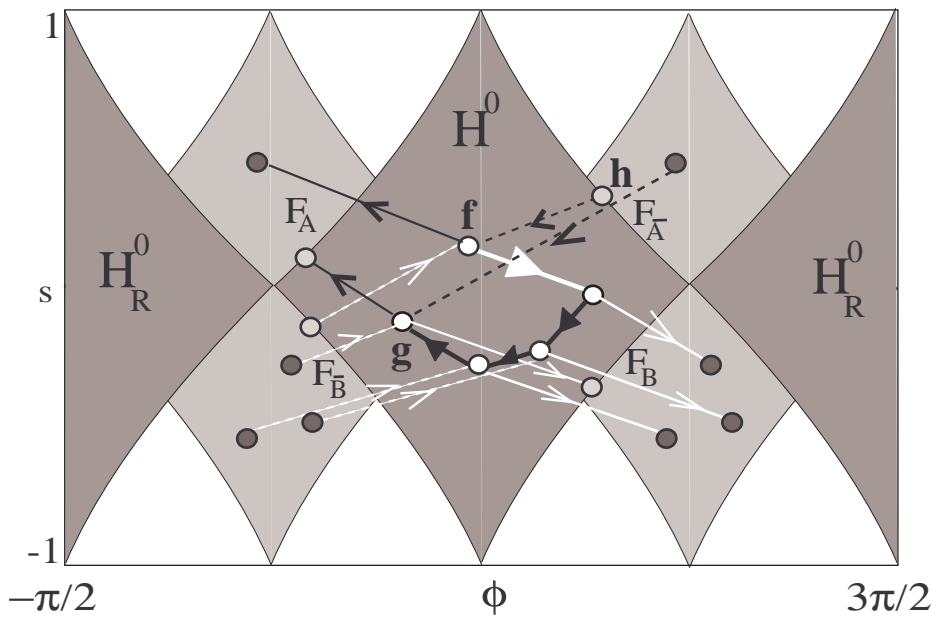}} \caption{Chain of points inside
$H^0$ and its terminating points on the frontier, for $\Delta=81$.
Black arrow: operator $A$, white arrow: operator $B$. The orbit is
asymmetric }\label{capp14}
\end{figure}

To conclude the proof of Theorem \ref{the6}, remember that, by
Lemma \ref{lem4}, all good points inside the interior of domains
$H_Z$  ($Z=A,B,\bar A, \bar B$) are image by $Z$ of points at the
interior of $H^0$. So, we are now interested in the images under
$Z^{-1}$ of the points of $F_Z$ which are inside $H^0$. For
instance, we start from a point ${\bf h}\in F_{\bar A}$ (see
Figure \ref{capp14}) and we go to ${\bf f}\in H^0$ applying  $A$.
Afterwards, we apply in sequel either $A$ or $B$ in order to
remain inside $H^0$ till we reach point ${\bf g}$,  such that both
$A{\bf g}$ and $B{\bf g}$ belong, by Lemma \ref{lem11}, to the
frontier of $H^0$, namely to $F_A$ and $F_B$. Lemma  \ref{lem11}
says also that ${\bar B}{\bf f}$ belongs to $F_{\bar B}$, since
${\bar A}{\bf f}$ belongs to $F_{\bar A}$. By a similar procedure
we associate to every point ${\bf h}$   of any one of sets $F_Z$,
a chain of points inside $H^0$ and  three other points of the
orbit of ${\bf f}$, one in everyone of the other sets $F_Y$,
$Y\not=Z$ (see Figure \ref{capp14}).  In his way  we   associate
to any point inside $H^0$ a unique ordered chain of points inside
$H^0$, whose initial point is mapped, by $\bar A$ and $\bar B$  to
two points of  $F_{\bar A}$ and $F_{\bar B}$, respectively, and
whose final point is mapped, by $A$ and $ B$ to two points of
$F_{A}$ and $F_{B}$, respectively (see figure \ref{capp14}). Since
also in this case different chains  cannot have common elements,
we have $\rho-1$ distinct orbits, corresponding to all integer
points in $F_{\bar A}$ ($(\rho,r,r)$ for  $r=1\dots\rho-1$), plus
the orbit of point $(\rho,0,0)$, given by Lemma \ref{lem12}.
\hfill $\square$

The above theorem has this immediate corollary:

\begin{cor} The good points in $H^0$ are partitioned into disjoint chains
whose points are obtained one from the preceding point by $A$ or
by $B$. The final point is sent by $A$ to $F_A$, by $B$ to $F_B$,
whereas the initial point is sent by $\bar A$ to $F_{\bar A}$ and
by $\bar B$ to $F_{\bar B}$.  Every chain corresponds to one
orbit.
\end{cor}

\begin{figure}[h]
\centerline{\epsfbox{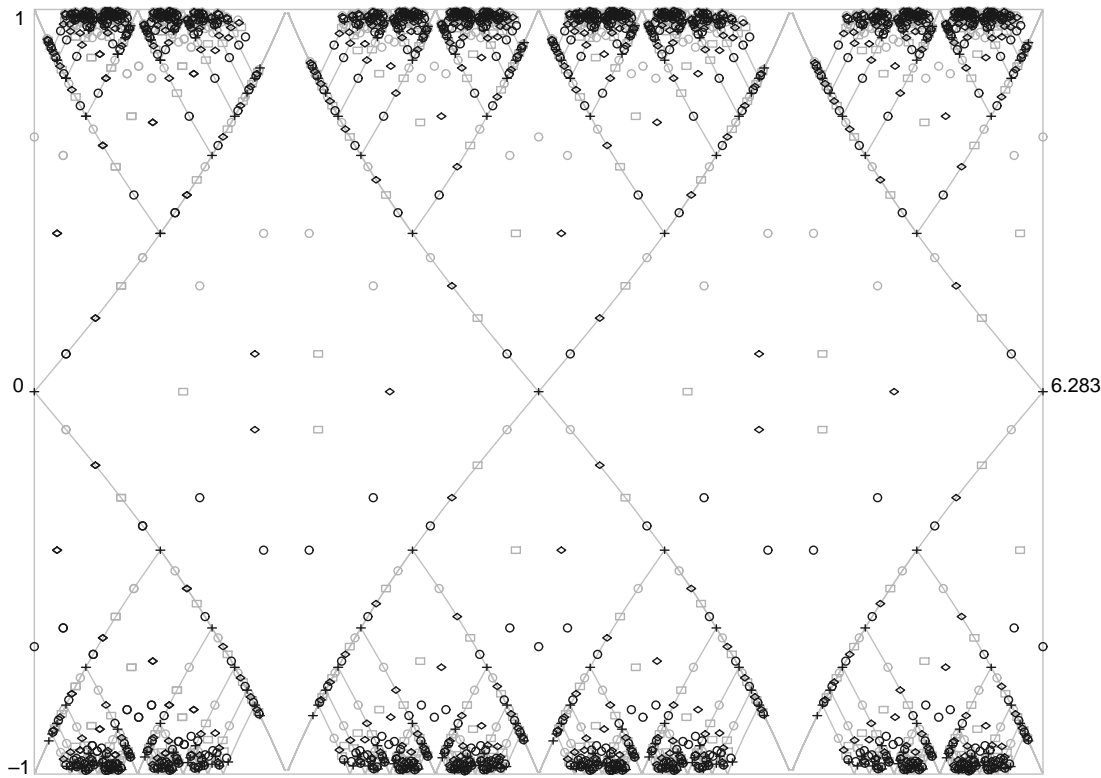}} \caption{ Initial part of the 5
distinct orbits for $\Delta=25$ projected on the cylinder. The
initial points of orbits are, in coordinates ($K,D,S$): (5,0,0),
crosses, supersymmetric orbit; (5,1,1), black circles,
$k$-symmetric orbit; (5,2,2), black diamonds, $(m+n)$-symmetric
orbit; (5,3,3), gray boxes, $(m+n)$-symmetric orbit; (5,4,4), gray
circles, $k$-symmetric orbit}\label{capp15}

\end{figure}

Figure \ref{capp15} shows the case $\Delta=25$.

\begin{thm}\label{the7} Let $t$ be the number of points of a chain
inside $H^0$, with initial point ${\bf f}$ and final point ${\bf
g}$. Let $T$ be the operator of $\mathcal{T}^+$ satisfying $T{\bf
f}={\bf g}$. Then $T$ is the product of $t-1$ generators. Among
them, $t_A$ are  of type $A$ and $t_B$ of type $B$, iff the orbit
of ${\bf f}$ contains exactly $t_B$ points inside every domain in
$G_A$ and $G_{\bar A}$, and exactly $t_A$ points inside every
domain in $G_B$ and $G_{\bar B}$, so that $t_A+t_B=t-1$.
\end{thm}

{\it Proof.} The points ${\bf f}={\bf f}_1,{\bf f}_2,\dots,{\bf
f}_{t}={\bf g}$ of the chain are mapped one to the successive one
by operator $A$ or $B$. By Lemma \ref{lem11}, to every image by
$A$ inside the chain there is an image by $B$ in $H_B$ and for
every image by $B$ inside the chain there is an image by $A$ in
$H_A$. Considering the reversal chain,  where point ${\bf f}$ is
reached from ${\bf g}$ by $t_A$ operators $\bar A$ and $t_B$
operators $\bar B$, we obtain that there are $t_A$ points  inside
$H_{\bar B}$ and $t_B$ points in $H_{\bar A}$. To complete the
proof we use Theorem \ref{the2}. \hfill $\square$

{\it Example.}  In Figure  \ref{capp15}, the orbit of (5,1,1)
(black circles) has a chain in $H_0$ composed of $4$ points. The
$3$ operators between them are all of type $A$, since there are no
black circles in $H_A$ and $H_{\bar A}$, and thus $t_A=3$ and
$t_B=0$.  Vice versa for the orbit of (5,4,4) (gray circles)
$t_A=0$ and $t_B=3$.  The orbits of (5,2,2) and (5,3,3), contain
in $H^0$ a chain of 3 points (diamonds and boxes, respectively),
and for both orbits $t_A=t_B=1$, since there is a diamond and a
box inside every domain, in regions $G,G_{\bar A}, G_B$ and
$G_{\bar B}$.

\section{Final remarks}

\subsection{Remark on the behaviour of the non rational points under $\mathcal T$.  }

Note that on the hyperboloid  $K^2+D^2-S^2=\Delta$ the orbit of
any point with fractional coordinates $(K,D,S)$, having the common
denominator $\mu$, is obtained (by a scale reduction) from the
orbit of the good point with integer coordinates $(2\mu K,2\mu
D,2\mu S)$ on the hyperboloid with discriminant $4\mu^2\Delta$,
and hence has a finite number of points in $H^0$, in $H^0_R$ as
well as in every domain of any generation.

Differently from  the elliptic case, where the orbit under
$\mathcal T$ of an irrational point is described exactly as that
of an integer point (close points in the Lobachevsky disc have
close orbits under ${\rm PSL}(2,\mathbb{Z})$ ), in the hyperbolic
case the situation is completely different.  Indeed, Theorems
\ref{the2} and \ref{the3} imply that two close points in $H$ have
close semi-orbits only  if these points  belong to the complement
of $H^0$ and $H_R^0$. But the orbits -- and even the semi-orbits
-- of two close points of $H^0$ or $H^0_R$ are not close (this
follows from the fact that an analogous statement as Lemma
\ref{lem6} holds for all points -- not only for the good points --
in $H^0$ and in $H_R^0$).

Moreover, for the irrational points on the one-sheeted hyperboloid
we obtain, from Theorems \ref{the2} and \ref{the3}, the following

\begin{cor} \label{cor3} The orbit of any point having at least one irrational
coordinate contains an infinite number of points in $H^0$ and in
$H^0_R$ (and hence in each connected component of the domains of
all generations of the hyperboloid). \end{cor}

\subsection{Sun eclipse model of the de Sitter world}$  \ \ \ \ \ $

In this section we see the Poincar\'e model of the de Sitter world
under an alternative projection.

Consider the domains of different generations directly on the
hyperboloid $H$. The lines bounding such domains belong to the
straight lines generatrices of the hyperboloid. More precisely, in
the plane $S=1$ the segment joining the point $p_i$ and its
opposite point on the circle $c_1$, upper vertices of a pair of
rhombi of the $n$-th generation, defines a direction on the plane
$S=1$. On the plane $S=0$, consider two straight lines $l_1$ and
$l_2$ in such direction tangent to the circle of radius
$\rho=\sqrt{\Delta}$, intersection of $H$ with that plane. The
four generatrices of the hyperboloid bounding  the domains
projected by $\mathcal Q$ to these rhombi are  the intersection of
the hyperboloid with two vertical planes through $l_1$ and $l_2$.

\begin{figure}[h]
\centerline{\epsfbox{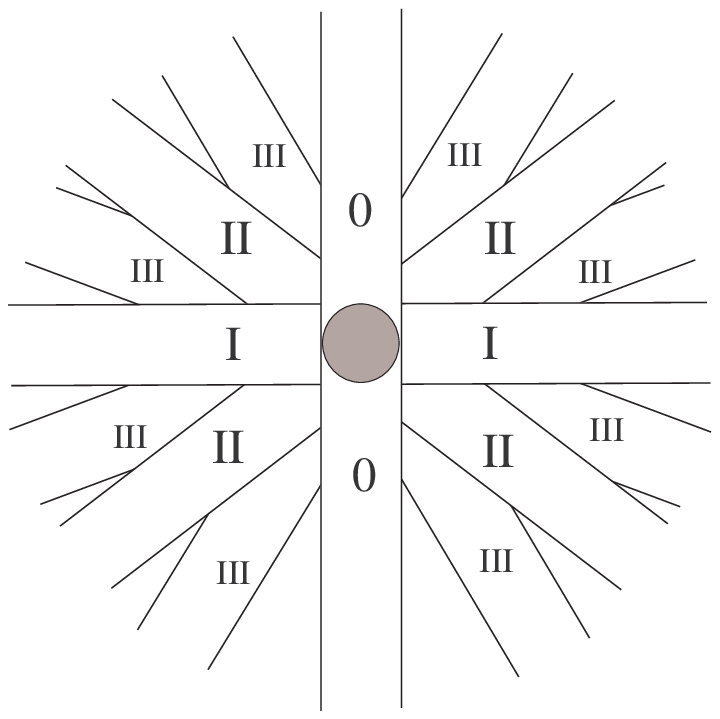}} \caption{}\label{capp16}
\end{figure}

The hierarchy of points $p_i$ is inherited by the pairs of
parallel lines on the plane $S=0$ as well as  by the regions
bounded by such pairs of lines and by the circle $K^2+D^2=\rho^2$.
The regions of the $n$-th generation {\it lie  behind} those of
all preceding generations. The view of the domains on the
hyperboloid projected to the plane $S=0$ is shown in Figure
\ref{capp16}. Note that by this projection  we map only one half
of $H$.

To introduce the `sun eclipse model', let us consider a property
of projection $\mathcal{P}$ (equation \ref{tilde}).

\begin{figure}[h]
\centerline{\epsfbox{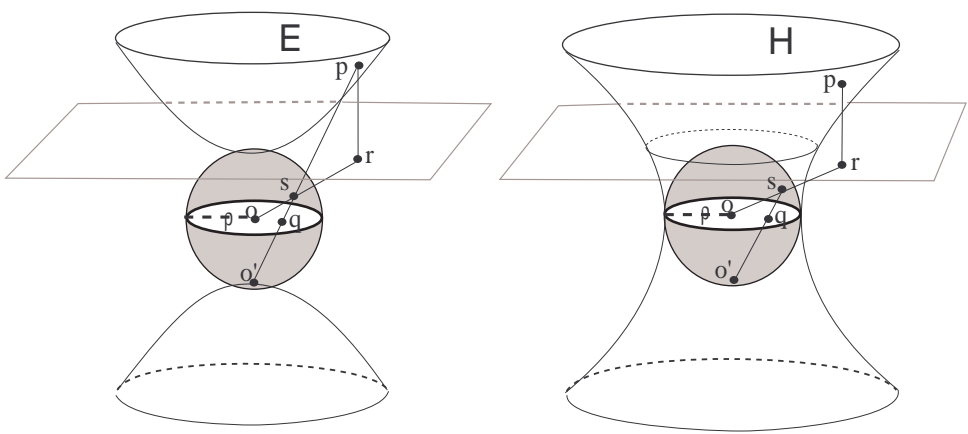}} \caption{Point $p$ on the
hyperboloid is sent to $q$ on the disc by a projection which
results  by the composition of three projections}\label{capp17}
\end{figure}

{\bf Definitions.} Let ${\bf f}$ be a point of the upper  sheet
$E$ of the two-sheeted hyperboloid ($K^2+D^2=S^2+\Delta$,
$\Delta<0$, see Figure \ref{capp17}). Let $\mathcal{P}'$ be the
projection of $E$ to the plane $S=\rho=\sqrt{-\Delta}$ along the
vertical direction, and ${\bf r}=\mathcal{P}' {\bf f} $.   Let
$\mathcal{P}''$ be the projection of the plane $S=\rho$ to the
sphere of radius $\rho$ from the centre of coordinates, and ${\bf
s}=\mathcal{P} ''{\bf r}$. Let $\mathcal{P}'''$ be the
stereographic projection of the upper half-sphere to the disc of
radius $\rho$ in the plane  $S=0$ from the point $O'$ ($K=D=0$,
$S=-\rho$), and ${\bf g}=\mathcal{P}''' {\bf s}$.

\begin{figure}[h]
\centerline{\epsfbox{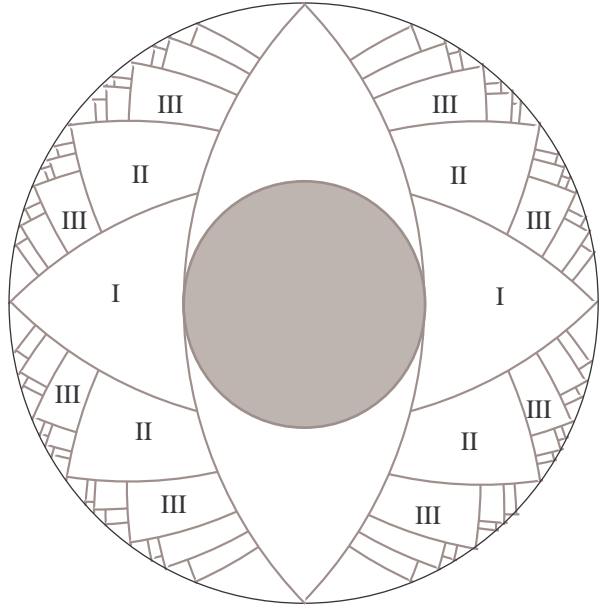}} \caption{Sun-eclipse model of the
de Sitter world }\label{capp18}
\end{figure}

\begin{prop}\label{pro3} Point $\bf q=\mathcal{P}'''\mathcal{P}''\mathcal {P}' {\bf p} $
coincides with the image of the projection  $\mathcal P{\bf p}$ of
${\bf p}$ directly to the plane $S=0$ from point $O'$. \end{prop}

{\it Proof.} Calculation. \hfill $\square$

The projection of the upper half-hyperboloid $H$
($K^2+D^2=S^2+\Delta$, $\Delta>0$) to the plane $S=0$ from point
$O'$ is not convenient at all: (it is two-to-one in a ring
contained in the unitary disc. However, projections
$\mathcal{P}'$, $\mathcal{P}''$ and $ \mathcal{P}'''$ (and their
symmetric ones for the lower half-hyperboloid) are well defined
(see Figure \ref{capp17}). We thus  project the upper
half-hyperboloid $H$ by $\mathcal{P}'''\mathcal{P}''\mathcal{P}'$
to the disc of radius $\rho$ to the plane $S=0$. The image is
contained in the ring $\frac {\rho}{2}<=\sqrt{K^2+D^2}< \rho$. The
pairs of straight lines tangent to the disc $K^2+D^2=\rho^2$ in
the plane $S=\rho$ are projected by $\mathcal{P}''$ to half
meridian circles of the sphere of radius $\rho$, and hence, by the
stereographic projection $\mathcal{P}'''$, to arcs of circles
tangent to the circle $\sqrt{k^2+D^2}=\rho/2$.   The final disc of
unit radius is obtained by  rescaling. On the boundary $C$ of this
disc the same points $p_i$ considered at the boundary of the
Lobachevsky disc are the extreme points of the domains of all
generations (forming the sun corona). The empty (black) disc of
radius $1/2$ is the moon in the sun-eclipse model (see Figure
\ref{capp18}).

{\it Remark.} The complementary forms (symmetric with respect the
$S$-axis) are in this model symmetric with respect to the centre
of the disc. Hence the picture of any orbit possesses this
symmetry. $k$-symmetric orbits  are symmetric with respect to the
vertical axis of the disc. A complete representation of an orbit,
either asymmetric or $k$-symmetric, requires two copies of this
model, one for the upper and the other for the lower half
hyperboloid.

\end{document}